\documentclass[10pt,letterpaper]{amsart}
\usepackage{amsmath, amssymb, amsthm, amsfonts}
\usepackage{mathtools}
\usepackage{graphicx}
\usepackage{subcaption}
\usepackage{fullpage}
\usepackage[dvipsnames]{xcolor}
\usepackage{tikz}
\usepackage{import,enumitem}
\usepackage{hyperref}
\usepackage[style=numeric, maxbibnames=15]{biblatex}
\renewbibmacro{in:}{}

\newtheorem{theorem}{Theorem}
\newtheorem{definition}[theorem]{Definition}
\newtheorem{lemma}[theorem]{Lemma}

\newtheorem{coro}[theorem]{Corollary}
\newtheorem{proposition}[theorem]{Proposition}
\newtheorem{remark}[theorem]{Remark}

\numberwithin{theorem}{section}
\theoremstyle{definition}
\numberwithin{equation}{section}
\usepackage{caption}
\usepackage{subcaption}
\newtheoremstyle{repeat}{}{}{\itshape}{}{\bfseries}{}{.5em}{\thmnote{#3}}
\theoremstyle{repeat}
\newtheorem*{repeattheorem}{Theorem}

\usetikzlibrary{patterns}
\newlength\cellsize 

\savebox2{
\begin{picture}(12,12)
\put(0,0){\line(1,0){12}}
\put(0,0){\line(0,1){12}}
\put(12,0){\line(0,1){12}}
\put(0,12){\line(1,0){12}}
\end{picture}}

\newcommand\cellify[1]{\def\thearg{#1}\def\nothing{}
\ifx\thearg\nothing\vrule width0pt height\cellsize depth0pt
  \else\hbox to 0pt{\usebox2\hss}\fi
  \vbox to 12\unitlength{\vss\hbox to 12\unitlength{\hss$#1$\hss}\vss}}

\newcommand\tableau[1]{\vtop{\let\\=\cr
\setlength\baselineskip{-12000pt}
\setlength\lineskiplimit{12000pt}
\setlength\lineskip{0pt}
\halign{&\cellify{##}\cr#1\crcr}}}

\def\cel(#1,#2)
{
  \node at (#1+0.5, #2+0.5) {$\times$};
  \draw (#1,#2) rectangle (#1+1, #2+1);
}
\def\celx(#1,#2)[#3]
{
  \node at (#1+0.5, #2+0.5) {$#3$};
  \draw (#1,#2) rectangle (#1+1, #2+1);
}

\def\rcel(#1,#2)
{
  \node[text=red] at (#1+0.5, #2+0.5) {$\times$};
  \draw[red] (#1,#2) rectangle (#1+1, #2+1);
}

\def\rcelx(#1,#2)[#3]
{
  \node[text=red] at (#1+0.5, #2+0.5) {$#3$};
  \draw[red] (#1,#2) rectangle (#1+1, #2+1);
}

\newcommand{\lock}{\reflectbox{$\mathbb{D}$}}

\newcommand{\N}{\mathbb{N}}
\newcommand{\Z}{\mathbb{Z}}

\DeclareMathOperator{\rwt}{rwt}
\DeclareMathOperator{\cwt}{cwt}

\DeclareMathOperator{\comp}{Comp}

\DeclareMathOperator{\movelabel}{move}

\definecolor{amethyst}{rgb}{0.6, 0.4, 0.8}

\title{Kohnert posets and polynomials of northeast diagrams}

\author[A.~Bingham]{Aram Bingham}\address[A.~Bingham]{Departmento de Matem\'aticas, Universidad de Chile, Ñuñoa, RM, Chile}
\email{\textcolor{blue}{\href{mailto:aram@matmor.unam.mx}{aram@matmor.unam.mx}}}

\author[B.~A.~Castellano]{Beth Anne Castellano}
\address[B.~A.~Castellano]{Department of Mathematics, Dartmouth College, Hanover, NH 03755 United States}
\email{\textcolor{blue}{\href{mailto:elizabeth.a.castellano.gr@dartmouth.edu}{elizabeth.a.castellano.gr@dartmouth.edu}}}

\author[K.~P.~Hadaway]{Kimberly P. Hadaway}
\address[K.~P.~Hadaway]{Department of Mathematics, Iowa State University, Ames, IA 50011 United States}
\email{\textcolor{blue}{\href{mailto:kph3@iastate.edu}{kph3@iastate.edu}}}

\author[R.~Hodges]{Reuven Hodges}
\address[R.~Hodges]{Department of Mathematics, University of Kansas, Lawrence, KS 66045 United States}
\email{\textcolor{blue}{\href{mailto:rmhodges@ku.edu}{rmhodges@ku.edu}}}

\author[Y.~Ma]{Yichen Ma}
\address[Y.~Ma]{Department of Mathematics, Texas State University, San Marcos, TX 78666 United States}
\email{\textcolor{blue}{\href{mailto:iwl24@txstate.edu}{iwl24@txstate.edu}}}

\author[A.~Moon]{Alex Moon}
\address[A.~Moon]{Department of Mathematics, Dartmouth College, Hanover, NH 03755 United States}
\email{\textcolor{blue}{\href{mailto:alexander.j.moon.gr@dartmouth.edu}{alexander.j.moon.gr@dartmouth.edu}}}

\author[K.~Salois]{Kyle Salois}
\address[K.~Salois]{Department of Mathematics, Statistics, and Computer Science, St. Olaf College, Northfield, MN 55057 United States}
\email{\textcolor{blue}{\href{mailto:salois1@stolaf.edu}{salois1@stolaf.edu}}}

\begin{document}

\maketitle{}

\begin{abstract}
Kohnert polynomials and their associated posets are combinatorial objects with deep geometric and representation theoretic connections, generalizing both Schubert polynomials and type-A Demazure characters. 
In this paper, we explore the properties of Kohnert polynomials indexed by \textbf{northeast diagrams}
along with their associated posets. 
We give separate classifications of the bounded, ranked, and multiplicity-free Kohnert posets for northeast diagrams, each of which can be computed in polynomial time with respect to the number of cells in the diagram. 
As an initial application, we specialize these classifications to simple criteria in the case of lock diagrams.
\end{abstract}

\section{Introduction} 

In \cite{AS22}, Assaf and Searles introduced \emph{Kohnert polynomials} which simultaneously generalize two fundamental families of polynomials, namely Schubert polynomials and (type-A) Demazure characters. Assaf and Searles define the Kohnert polynomials as generating polynomials of certain diagrams in the plane. They drew inspiration from the work of Kohnert in \cite{K91} that gives such a generating polynomial interpretation for the Demazure characters. In his thesis, Kohnert also conjectured that the Schubert polynomials had such an interpretation, and this was eventually proven in \cite{A22,W02,W99}. 

Schubert polynomials have their origin in the geometry of flag varieties. First introduced by Lascoux and Sch\"{u}tzenberger, Schubert polynomials are polynomial representatives of Schubert classes in the cohomology ring of type-A flag varieties. It is a central problem in algebraic combinatorics to give a combinatorial interpretation of the structure constants for multiplication of Schubert polynomials~\cite{P24, S00}. On the other hand, Demazure characters (also known as key polynomials) have connections to both representation theory and geometry. Demazure characters are the characters of certain finite dimensional representations, known as Demazure modules, that arise as the space of sections of a line bundle on the type-A flag variety that has been restricted to a Schubert variety. Noting these deep connections to geometry and representation theory, Assaf gave a construction in \cite{A21} of a Demazure crystal whose character is the Kohnert polynomial whenever the indexing diagram is \textbf{southwest}.

Given the discussion of the preceding paragraph, it should come as no surprise that the coefficients of Schubert polynomials and Demazure characters in the monomial basis have geometric implications. If all such coefficients are equal to $0$ or $1$ then we say that these polynomials are \textbf{monomial multiplicity-free}. In \cite{FMS21}, Fink, M\'{e}sz\'{a}ros, and St. Dizier gave a classification of the monomial multiplicity-free Schubert polynomials in terms of permutation pattern avoidance. 
In \cite{HY23}, the fourth author and Yong gave a classification of the monomial multiplicity-free Demazure characters in terms of composition pattern avoidance. The monomial multiplicity-freeness of these polynomials is equivalent, via character theory, to the associated algebraic torus module being multiplicity-free as a direct sum of irreducible torus modules. This in turn has implications for the study of torus orbits in Schubert varieties~\cite{HY22}. Using these multiplicity-freeness results Gao, the fourth author, and Yong gave a classification of the Schubert varieties in the flag variety that are spherical varieties for the action of a Levi subgroup of the general linear group~\cite{GHY23}. 

With the wide ranging connections to fundamental problems exhibited by these Kohnert polynomials and their special cases, a better understanding of the defining combinatorial model is imperative. To discuss results in this direction, we now explicate the Kohnert polynomial construction in detail. 

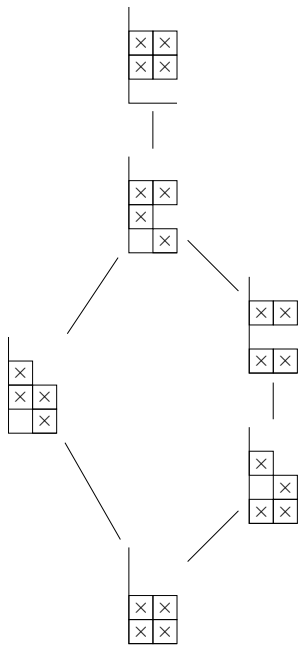
\begin{figure}[ht]
    \centering
    \scalebox{0.8}{\begin{tikzpicture}
 \node (1) at (0,0) {\begin{tikzpicture}[scale=0.4]
    \cel(0,0)
    \cel(1,0)
    \cel(0,1)
    \cel(1,1)
    \draw (0,4)--(0,0)--(2,0);
\end{tikzpicture}};
\node (2) at (-2,3.5) {\begin{tikzpicture}[scale=0.4]
   \cel(0,2)
    \cel(1,0)
    \cel(0,1)
    \cel(1,1)
  \draw (0,4)--(0,0)--(2,0);
\end{tikzpicture}};
\node (3) at (2,2) {\begin{tikzpicture}[scale=0.4]
  \cel(0,0)
    \cel(1,0)
    \cel(0,2)
    \cel(1,1)
  \draw (0,4)--(0,0)--(2,0);
\end{tikzpicture}};
\node (4) at (2,4.5) {\begin{tikzpicture}[scale=0.4]
     \cel(0,0)
    \cel(1,0)
    \cel(0,2)
    \cel(1,2)
  \draw (0,4)--(0,0)--(2,0);
\end{tikzpicture}};
\node (5) at (0,6.5) {\begin{tikzpicture}[scale=0.4]
    \cel(0,2)
    \cel(1,0)
    \cel(0,1)
    \cel(1,2)
  \draw (0,4)--(0,0)--(2,0);
\end{tikzpicture}};
\node (6) at (0,9) {\begin{tikzpicture}[scale=0.4]
 \cel(0,2)
    \cel(1,2)
    \cel(0,1)
    \cel(1,1)
  \draw (0,4)--(0,0)--(2,0);
\end{tikzpicture}};
\draw (2)--(1)--(3);
\draw (2)--(5);
\draw (3)--(4)--(5);
\draw (5)--(6);
\end{tikzpicture}}
    \caption{$\mathcal{P}(D)$ where $D=\{(2,1),(2,2),(3,1),(3,2)\}$ }\label{fig:basicposet}
\end{figure}

A \textbf{diagram} $D$ is a finite subset of $\N\times \N$. We refer to $(r,c) \in D$ as a cell with row and column index $r$ and $c$, respectively. Equivalently, we visualize a diagram as a collection of row and column coordinates of cells with rows labeled from bottom to top and columns labeled left to right (see Figure \ref{fig:basicdiagram} below).
\begin{figure}[ht]
\begin{center}    
    \begin{tikzpicture}[scale=0.5]
    \cel(1,1)
    \cel(0,2)
    \cel(0,3)
    \cel(1,3)
    \draw (0,5)--(0,0)--(2,0);
    \end{tikzpicture}
   \caption{The diagram $D=\{(2,2),(3,1),(4,1),(4,2)\}$}\label{fig:basicdiagram} 
  \end{center}
\end{figure}
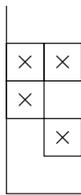
A \textbf{Kohnert move} in a diagram $D$ selects the rightmost cell in a row and moves that cell down to the first empty position below it, jumping over other cells if needed. A Kohnert move is called \textbf{elementary} if the cell moves to the position directly below itself, not jumping over any other cells. Otherwise, it is called a \textbf{jump} Kohnert move. The \textbf{Kohnert poset} of a diagram $D$, denoted by $\mathcal{P}(D)$, has an underlying set consisting of all diagrams that can be obtained from $D$ by the application of a (possibly empty) sequence of Kohnert moves. It is the transitive closure of the relations $D_2 < D_1$ for $D_1, D_2 \in \mathcal{P}(D)$ if $D_2$ is the result of applying a single Kohnert move to $D_1$. Given two diagrams $D$ and $E$, we write $D \leq E$ to indicate that $D \in \mathcal{P}(E)$. See Figure \ref{fig:basicposet} for an example Kohnert poset.

We remark that since Kohnert moves only impact individual columns, empty columns in a diagram have no bearing on its Kohnert poset structure. Therefore, without loss of generality, in this paper we assume all diagrams have no empty columns. 

Towards associating polynomials with diagrams, let $\comp_n := \mathbb{Z}^{n}_{\geq 0}$ be the set of weak compositions of length $n$. The \textbf{row weight} of a diagram $D$, denoted by $\rwt(D)$, is the weak composition $(\alpha_1,\ldots,\alpha_{\ell})$ such that $\alpha_i$ is the number of cells in row $i$, and $\ell$ is the maximum row index such that row $\ell$ is nonempty.
The \textbf{column weight} $\cwt(D)$ is the composition $(\beta_1,\ldots, \beta_m)$ such that $\beta_i$ is the number of cells in column $i$, and column $m$ is the rightmost nonempty column in $D$. In the example in Figure~\ref{fig:basicdiagram}, $\rwt(D) = (0,1,1,2)$ and $\cwt(D) = (2,2)$.

\begin{definition}[{\cite[Definition 2.2]{AS22}}] 
Given a diagram $D$, its \textbf{Kohnert polynomial} is
\begin{center}
$\displaystyle \mathfrak{K}_D = \sum_{T \in \mathcal{P}(D)} x_1^{\rwt(T)_1} \cdots x_n^{\rwt(T)_n}$. 
\end{center}
\end{definition}

In \cite{AS22}, Assaf and Searles show that if $\{ D_{\alpha} \}$ is a set of diagrams indexed by weak compositions with the property that $\rwt(D_{\alpha})=\alpha$ for $\alpha \in \comp_n$, then their Kohnert polynomials $\{ \mathfrak{K}_{D_{\alpha}} \}$ form a basis of the polynomial ring $\mathbb{Z}[x_1,\ldots,x_n]$. As remarked earlier, these bases subsume both the Schubert and Demazure character bases of the polynomial ring.

Now, moving towards the content of the present work, Colmenarejo,  Hutchins, Mayers, and Phillips initiated a comprehensive study of the structural properties of $\mathcal{P}(D)$ in \cite{CHMP23}. Therein, they determine a sufficient condition guaranteeing when such Kohnert posets are bounded and
two necessary conditions for when they are ranked. They also show that the Kohnert posets associated with Demazure characters are always bounded and give a criterion for when the poset is ranked. Inspired by their work, we investigate the Kohnert posets of \textbf{northeast diagrams} (see Definition~\ref{def:northeast}). Our initial motivation was to investigate a subclass of northeast diagrams called \textbf{lock diagrams}. Lock diagrams were first defined by Assaf and Searles~\cite{AS22} as natural analogs of the left-justified diagrams of Demazure characters. The polynomials associated to lock diagrams are called \textbf{lock polynomials}. In \cite{W20}, Wang showed that lock polynomials exhibit a Demazure crystal structure that is closely intertwined with the crystal structure of Demazure characters. 

Here, we give a complete classification of the northeast diagrams $D$ for which $\mathcal{P}(D)$ is ranked and, separately, bounded. Additionally, we give a complete criterion for when $\mathfrak{K}_D$ is monomial multiplicity-free for any northeast diagram $D$. A highlight of all three classifications is that they are given purely in terms of the diagram $D$; that is, they do not require any analysis of the underlying set of diagrams in $\mathcal{P}(D)$.  Along the way, we develop a suite of tools that aid in the study of Kohnert posets in more general settings.

The feasibility of all three of our classifications is based on the following surprising result, which comprises our first main theorem. The underlying set of the \textbf{elementary Kohnert poset} $\mathcal{P}^{ele}(D)$ consists of all diagrams that can be obtained from $D$ by the application of a (possibly empty) sequence of elementary Kohnert moves. It is the transitive closure of the relations $D_2 < D_1$ for $D_1, D_2 \in \mathcal{P}(D)$ if $D_2$ is the result of applying a single \textbf{elementary} Kohnert move to $D_1$.

\begin{theorem} \label{thm:NEKohnertPosetEle}
Let $D$ be a northeast diagram. Then $\mathcal{P}(D)$ is a refinement of $\mathcal{P}^{ele}(D)$; in particular, $\mathcal{P}(D)$ and $\mathcal{P}^{ele}(D)$ are equal as sets.
\end{theorem}

In other words, any diagram $D_1 \in \mathcal{P}(D)$ may be reached from $D$ via only elementary moves. This non-trivial fact makes the analysis of the Kohnert posets of northeast diagrams considerably more tractable. The converse of Theorem~\ref{thm:NEKohnertPosetEle} is false in general. For instance, any diagram with at most one cell in each column satisfies that $\mathcal{P}(D)$ is a refinement of $\mathcal{P}^{ele}(D)$, regardless of whether $D$ is northeast. It would be very instructive to classify the diagrams $D$ for which this phenomenon occurs. In Section~\ref{sec:labelings}, we develop the notion of tableaux, expanding on the idea of Kohnert tableaux introduced by Assaf in \cite[Section 5]{A21}, in a setting more general than northeast diagrams. These tableaux are instrumental in proving the subsequent theorems.

We now state our classification results. First, our classification of the monomial multiplicity-free Kohnert polynomials of northeast diagrams:
\begin{theorem}\label{thm:mainMFree}
    If $D_0$ is a northeast diagram, then $\mathfrak{K}_{D_0}$ is monomial multiplicity-free if and only if $D_0$ does not contain $x_1 = (r_1,c_1)$ and $x_2 = (r_2,c_2)$ such that:
    \begin{itemize}
        \item[(a)] $r_1 < r_2$,
        \item[(b)] $c_1 < c_2$,
        \item[(c)] there is at least one empty position $(r, c_1)$ where $r < r_1$, and
        \item[(d)] for each $c>c_1$, there are at least two empty positions $(r,c)$ where $r \leq r_1$. 
    \end{itemize}
\end{theorem}

We say that a poset $P$ is \textbf{ranked} if there exists a rank function $\rho: P \rightarrow \mathbb{Z}$ such that if $y$ covers $x$, then $\rho(y) = \rho(x) + 1$. This definition does not require all minimal elements to have rank $0$. An equivalent characterization of the poset $P$ being ranked is that for all $x \le y$ in $P$, every maximal chain from $x$ to $y$ has the same length.

\begin{theorem}\label{thm:mainNEranked}
     If $D_0$ is a northeast diagram, then $\mathcal{P}(D_0)$ is ranked if and only if $D_0$ does not contain $x_1 = (r_1,c_1), x_2 = (r_2,c_2)$, and $x_3 = (r_3,c_3)$ such that:
 \begin{itemize}
         \item[(a)] $r_1 < r_2 \leq r_3$,
         \item[(b)] $c_1 = c_2 < c_3$, 
         \item[(c)] for each  $c_1 \leq c <c_3$, there is at least one empty position $(r,c)$ where $r < r_1$, and
         \item[(d)] for each  $c \geq c_3$, the number of $r < r_3$ such that $(r, c) \in D_0$ is strictly less than $r_1$.
        \end{itemize}
\end{theorem}
\noindent 

\noindent A poset $P$ is \textbf{bounded} if $P$ contains a unique minimal element and a unique maximal element. The Kohnert poset $\mathcal{P}(D)$ always has a unique maximal element, namely $D$. Hence, a Kohnert poset is bounded if and only if it has a unique minimal element.

\begin{theorem}\label{thm:mainNEbounded}
    If $D_0$ is a northeast diagram, then $\mathcal{P}(D_0)$ is bounded if and only if $D_0$ does not contain $x_1 = (r_1,c_1), x_2 = (r_2,c_2)$, and $x_3 = (r_3,c_3)$  such that:
    \begin{itemize}
        \item[(a)] $r_1 \leq r_2 < r_3$,
        \item[(b)] $c_1 < c_2 = c_3$,
        \item[(c)] for each $c_1 \leq c < c_2$, $\cwt(D_0)_c < \cwt(D_0)_{c_2}$,
        \item[(d)] for each  $c\geq c_1$, there is at least one empty position $(r,c)$ where $r < r_1$, and 
        \item[(e)] for each $r_1 < r \leq r_3$, the cell $(r,c_1)$ is not in $D_0$.
    \end{itemize}
\end{theorem}

As an initial application of these results we specialize each of these classifications to lock diagrams and polynomials, yielding significantly simpler criterion. In future work, we plan to focus on lock diagrams in greater detail and study shellability, rank-unimodality, and enumerative formulas for the number of minimal elements in their Kohnert posets. 

The organization of this paper is as follows. 
In Section \ref{sec:tools}, we define northeast diagrams, introduce labelings on northeast diagrams and the diagrams that arise in their posets, and prove several technical results that are used in the arguments for our main theorems. 
In Sections \ref{sec:Mfree}, \ref{sec:ranked}, and \ref{sec:bounded}, we give the proofs of Theorems \ref{thm:mainMFree}, \ref{thm:mainNEranked}, and \ref{thm:mainNEbounded}, respectively. Of the three, we note that the proof of Theorem \ref{thm:mainNEbounded} is the most complex.
Finally, in Section \ref{sec:lock} we describe and prove how these results specialize to \emph{lock diagrams}.

\section{Tools for diagrams}\label{sec:tools}
The purpose of this section is to define our diagrams of interest and develop the tools required to prove our three main classification theorems. Some of these notions and results will also be useful in our forthcoming work that studies combinatorial properties of the posets associated with lock diagrams in greater detail, and they may be of general interest to researchers working on Kohnert diagram combinatorics.

\subsection{Classes of diagrams}\label{sec:diagrams}

We begin by defining the class of diagrams that form our primary objects of study.

\begin{definition}\label{def:northeast}
    A diagram $D$ is \textbf{northeast} if for all pairs of cells $(r_1, c_1), (r_2, c_2) \in D$, we have that $(\max(r_1, r_2), \max(c_1, c_2)) \in D$. Similarly, $D$ is \textbf{southeast} if for all pairs of cells $(r_1, c_1), (r_2, c_2) \in D$, we have that $(\min(r_1, r_2), \max(c_1, c_2)) \in D$.
\end{definition}

\begin{figure}[ht]
\begin{subfigure}{0.45\textwidth}
  \begin{center}
    \begin{tikzpicture}[scale=0.5]
    \cel(0,1)
    \cel(1,2)
    \cel(0,3)
    \cel(1,3)
    \cel(1,4)
    \cel(2,4)
    \draw (0,6)--(0,0)--(3,0);
    \end{tikzpicture}
    \caption{ A northeast diagram}
  \end{center}
\end{subfigure}
\begin{subfigure}{0.45\textwidth}
  \begin{center}
    \begin{tikzpicture}[scale=0.5]
    \cel(1,1)
    \cel(2,1)
    \cel(0,2)
    \cel(1,2)
    \cel(2,2)
    \cel(0,3)
    \cel(1,3)
    \cel(0,4)
    \draw (0,6)--(0,0)--(3,0);
    \end{tikzpicture}
    \caption{ A southeast diagram}
  \end{center}
\end{subfigure} \caption{A northeast (left) and a southeast (right) diagram.}\label{fig:NESEex}
\end{figure}

See Figure \ref{fig:NESEex} for examples. The naming of these diagrams arises from intuition relating to the cardinal directions.
In a northeast diagram, for every two cells where one cell (in position $(r_1,c_1)$) is strictly northwest of the other (in position $(r_2,c_2)$), there is a cell located at the ``northeast corner"  of the rectangle defined by the two cells, i.e., the position $(r_1,c_2)$. Similarly, in a southeast diagram, for every two cells where one cell is strictly southwest of the other, there is a cell located at the ``southeast corner" of the rectangle defined by the two cells. It is possible to define northwest and southwest diagrams in a similar manner, but we omit the definitions since they are not needed in this work. 

An important family of diagrams is \textbf{lock diagrams}, which are exactly those diagrams that are both northeast and southeast.

\begin{definition}\label{def:lock_diagram}
    For $\alpha = (\alpha_1, \ldots, \alpha_n) \in \comp_n$ with $m = \max (\alpha_1, \dots, \alpha_n)$, the \textbf{lock diagram} is the diagram
    \[\lock(\alpha) \coloneqq \{(r, c) \mid 1 \leq r \leq n, m - \alpha_r \leq c \leq m\}.\]
    $\lock(\alpha)$ is said to be the \textbf{lock diagram of row weight $\alpha$.}
\end{definition}

\begin{figure}[ht]
  \begin{center}
    \begin{tikzpicture}[scale=0.5]
    \cel(2,0) 
    \cel(0,1) 
    \cel(1,1) 
    \cel(2,1)
    \cel(2,2)
    \cel(1,4)
    \cel(2,4)
    \draw (0,5)--(0,0)--(3,0);
    \end{tikzpicture}
    \caption{ The lock diagram $\lock(\alpha)$ for $\alpha = (1, 3, 1, 0 , 2)$ .}\label{fig:lockdiagram}
  \end{center}
\end{figure}
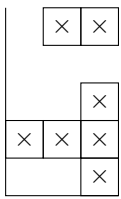

Eqvuialently, $\lock(\alpha)$ has exactly $\alpha_r$ cells in row $r$, and if $c < m$ with $(r, c) \in \lock(\alpha)$, then $(r, c+1) \in \alpha$ as well. See Figure \ref{fig:lockdiagram} for an example of a lock diagram. 

\begin{lemma}\label{lemma:northeastSoutheast}
    A diagram $D$ is both northeast and southeast if and only if it is a lock diagram.
\end{lemma}

\begin{proof}
Assume that $D$ has $m$ columns.

        (Proof of $\Rightarrow$)  Suppose $D$ satisfies both the northeast and southeast conditions and that $(r, c) \in D$ for some $1 \leq c < m$. As column $c+1$ is nonempty by assumption, there exists some $r'$ such that $(r', c+1) \in D$. We have three cases.

    \begin{itemize}[leftmargin=0.65in]
\item[\textbf{Case 1:}] $r' > r$. Then $(r, c+1) \in D$ by the southeast condition.
    
   \item[\textbf{Case 2:}] $r' = r$. Then $(r, c+1) \in D$ by assumption.

    \item[\textbf{Case 3:}] $r' < r$. Then $(r, c+1) \in D$ by the northeast condition.
    \end{itemize}
    Therefore, by definition, $D$ is a lock diagram.
    
    (Proof of $\Leftarrow$)  Suppose $D$ is a lock diagram and $(r_1, c_1), (r_2, c_2) \in D$. By definition of a lock diagram, $(r_1, \max(c_1,c_2))$ and $(r_2, \max(c_1,c_2))$ are both in $D$. Thus, $D$ satisfies both the northeast and southeast conditions.
\end{proof}

\subsection{Labelings and tableaux}\label{sec:labelings}

In what follows, we develop the notions of northeast labelings and tableaux. This culminates with Theorem \ref{thm:labeling}, which yields a tableau criterion for a diagram to be contained in the Kohnert poset of a northeast diagram.

\begin{definition}
A \textbf{labeling} of a diagram $D$ is a map $\mathcal{L}: D \to \N$. A labeling is \textbf{strict} if the labels within each column are strictly increasing from bottom to top. The \textbf{super-standard} labeling of a diagram labels each cell with its row index. 
\end{definition}

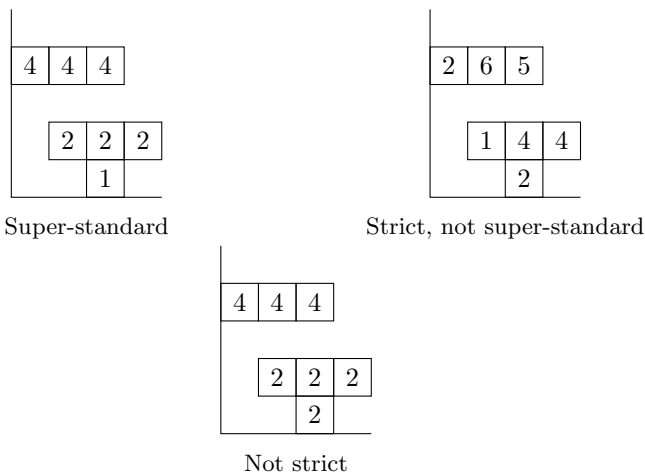
\begin{figure}[ht]
\begin{subfigure}{0.33\textwidth} \centering
    \begin{tikzpicture}[scale=0.5]
    \draw (2,0) rectangle (3,1); \node at (2.5, 0.5) {$1$};%
    \draw (1,1) rectangle (2,2); \node at (1.5, 1.5) {$2$};%
    \draw (2,1) rectangle (3,2); \node at (2.5, 1.5) {$2$};%
    \draw (3,1) rectangle (4,2); \node at (3.5, 1.5) {$2$};%
    \draw (0,3) rectangle (1,4); \node at (0.5, 3.5) {$4$};%
    \draw (1,3) rectangle (2,4); \node at (1.5, 3.5) {$4$};%
    \draw (2,3) rectangle (3,4); \node at (2.5, 3.5) {$4$};%
    \draw (0,5)--(0,0)--(4,0);
    \end{tikzpicture}
    \caption{Super-standard
    }\label{subfig:superstandard}
\end{subfigure}
\begin{subfigure}{0.33\textwidth}\centering
    \begin{tikzpicture}[scale=0.5]
    \draw (2,0) rectangle (3,1); \node at (2.5, 0.5) {$2$};%
    \draw (1,1) rectangle (2,2); \node at (1.5, 1.5) {$1$};%
    \draw (2,1) rectangle (3,2); \node at (2.5, 1.5) {$4$};%
    \draw (3,1) rectangle (4,2); \node at (3.5, 1.5) {$4$};%
    \draw (0,3) rectangle (1,4); \node at (0.5, 3.5) {$2$};%
    \draw (1,3) rectangle (2,4); \node at (1.5, 3.5) {$6$};%
    \draw (2,3) rectangle (3,4); \node at (2.5, 3.5) {$5$};%
    \draw (0,5)--(0,0)--(4,0);
    \end{tikzpicture}
    \caption{Strict, not super-standard}\label{subfig:strict}
\end{subfigure}
\begin{subfigure}{0.33\textwidth} \centering
    \begin{tikzpicture}[scale=0.5]
    \draw (2,0) rectangle (3,1); \node at (2.5, 0.5) {$2$};%
    \draw (1,1) rectangle (2,2); \node at (1.5, 1.5) {$2$};%
    \draw (2,1) rectangle (3,2); \node at (2.5, 1.5) {$2$};%
    \draw (3,1) rectangle (4,2); \node at (3.5, 1.5) {$2$};%
    \draw (0,3) rectangle (1,4); \node at (0.5, 3.5) {$4$};%
    \draw (1,3) rectangle (2,4); \node at (1.5, 3.5) {$4$};%
    \draw (2,3) rectangle (3,4); \node at (2.5, 3.5) {$4$};%
    \draw (0,5)--(0,0)--(4,0);
    \end{tikzpicture}
    \caption{Not strict}\label{subfig:notstrict}
\end{subfigure}
\caption{Three diagrams with labelings.}\label{fig:labelings}
\end{figure}

See Figure \ref{fig:labelings} for examples of these types of labelings. We refer to the pair of a diagram and its labeling $T= (D,\mathcal{L})$ as a \textbf{tableau} (this differs from the tableaux of \cite{AS18}).

\begin{definition} 
    The \textbf{column content} $C_\mathcal{L} = (C_1, C_2, \ldots)$ of a strict labeling $\mathcal{L}: D \to \mathbb{N}$ is the sequence such that $C_i \subseteq \mathbb{Z}_{\geq 0}$ is set of labels of cells in column $i$. Two labelings $\mathcal{L}: D \to \mathbb{N}$ and $\mathcal{L}': D' \to \mathbb{N}$ are said to be \textbf{column-equivalent} if $C_\mathcal{L} = C_{\mathcal{L}'}$. We also refer to two tableaux with column-equivalent labelings as column-equivalent in their own right.
\end{definition}

In other words, two labelings are column-equivalent if they have the same set of labels in each fixed column.

\begin{definition}
    Given a diagram $D \in \mathcal{P}(D_0)$, the \textbf{standard labeling $\mathcal{L}$ of $D$ with respect to $D_0$} is the unique strict labeling of $D$ that is column-equivalent to the super-standard labeling of $D_0$. When the poset containing $D$ is clear from context, we call $\mathcal{L}$ the standard labeling of $D$. In a diagram with a standard labeling, we refer to a cell with the label $i$ as an \textbf{$i$-cell.}
\end{definition}

\begin{lemma} \label{lem:stdLabelingsConstruction} 
If $D_b<D_a$ in $\mathcal{P}(D_0)$, then the standard labeling of $D_b$ with respect to $D_0$ can be constructed from the standard labeling of $D_a$ with respect to $D_0$. Moreover, if $D_b$ can be obtained from $D_a$ by a single Kohnert move, then for any $x = (r, c) \in D_a$, the cell $x' = (r', c) \in D_b$ with the same label as $x$ has $r - 1 \leq r' \leq r$.
\end{lemma}

\begin{proof}
In the case that $D_b$ results from $D_a$ by a single Kohnert move, then the standard labeling of $D_b$ can be constructed as follows. The diagrams $D_a$ and $D_b$ differ by the placement of exactly one cell which we denote as $x_1$. We address two possible cases. If the Kohnert move that yields $D_b$ from $D_a$ was elementary (and therefore consisted of the cell $x_1$ moving one position down into an empty space), then preserve the label on $x_1$. If the Kohnert move was a jump move, then $x_1$ jumped over a string of adjacent cells $x_2,x_3,\ldots, x_{\ell}$, where cell $x_i$ is directly above $x_{i+1}$ for all $i$. In this case, we relabel $x_2$ with the label for $x_1$, relabel $x_3$ with the label for $x_2$, and so on, finally relabeling $x_1$ with the label for $x_{\ell}$. The resulting labeling is the standard labeling of $D_b$. In this sense, we re-frame a Kohnert move as ``pushing'' all of the cells directly below it into the next available empty space, preserving the labels on all moved cells.

On the other hand, if $D_b$ and $D_a$ differ by more than one Kohnert move, then the standard labeling of $D_b$ may be constructed from the standard labeling $\mathcal{L}_a$ of $D_a$ by taking any chain $D_b=C_{1}<C_{2}<\cdots<C_{k}=D_a$ such that $C_k$ results from $C_{k+1}$ by a single Kohnert move for all $k$. Then one can repeatedly apply the process described in the preceding paragraph, starting with $\mathcal{L}_a$. The choice of chain does not affect this process since standard labelings are unique. Further, the super-standard labeling of $D_0$ is the standard labeling of $D_0$, and so this process gives a way of constructing the standard labelings of all diagrams in $\mathcal{P}(D_0)$. 
\end{proof}

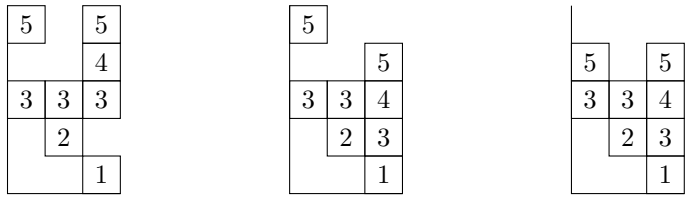
\begin{figure}[ht] \begin{center}
\begin{tikzpicture}[scale = 0.5]
    \draw (0,5)--(0,0)--(3,0);
    \draw (2,0) rectangle (3,1); \node at (2.5, 0.5) {$1$};%
    \draw (1,1) rectangle (2,2); \node at (1.5, 1.5) {$2$};%
    \draw (0,2) rectangle (1,3); \node at (0.5, 2.5) {$3$};%
    \draw (1,2) rectangle (2,3); \node at (1.5, 2.5) {$3$};%
    \draw (2,2) rectangle (3,3); \node at (2.5, 2.5) {$3$};%
    \draw (2,3) rectangle (3,4); \node at (2.5, 3.5) {$4$};%
    \draw (0,4) rectangle (1,5); \node at (0.5, 4.5) {$5$};%
    \draw (2,4) rectangle (3,5); \node at (2.5, 4.5) {$5$};%
\end{tikzpicture} \hspace{2cm}
\begin{tikzpicture}[scale = 0.5]
    \draw (0,5)--(0,0)--(3,0);
    \draw (2,0) rectangle (3,1); \node at (2.5, 0.5) {$1$};%
    \draw (1,1) rectangle (2,2); \node at (1.5, 1.5) {$2$};%
    \draw (0,2) rectangle (1,3); \node at (0.5, 2.5) {$3$};%
    \draw (1,2) rectangle (2,3); \node at (1.5, 2.5) {$3$};%
    \draw (2,1) rectangle (3,2); \node at (2.5, 1.5) {$3$};%
    \draw (2,2) rectangle (3,3); \node at (2.5, 2.5) {$4$};%
    \draw (0,4) rectangle (1,5); \node at (0.5, 4.5) {$5$};%
    \draw (2,3) rectangle (3,4); \node at (2.5, 3.5) {$5$};%
\end{tikzpicture} \hspace{2cm}
\begin{tikzpicture}[scale = 0.5]
    \draw (0,5)--(0,0)--(3,0);
    \draw (2,0) rectangle (3,1); \node at (2.5, 0.5) {$1$};%
    \draw (1,1) rectangle (2,2); \node at (1.5, 1.5) {$2$};%
    \draw (0,2) rectangle (1,3); \node at (0.5, 2.5) {$3$};%
    \draw (1,2) rectangle (2,3); \node at (1.5, 2.5) {$3$};%
    \draw (2,1) rectangle (3,2); \node at (2.5, 1.5) {$3$};%
    \draw (2,2) rectangle (3,3); \node at (2.5, 2.5) {$4$};%
    \draw (0,3) rectangle (1,4); \node at (0.5, 3.5) {$5$};%
    \draw (2,3) rectangle (3,4); \node at (2.5, 3.5) {$5$};%
\end{tikzpicture}
\caption{The standard labelings of two diagrams (center, right) with respect to a diagram with a super-standard labeling (left).}\label{fig:standard_label_ex} \end{center}
\end{figure}

We reference Figure \ref{fig:standard_label_ex} for an example of this process. 

In the setting of labelings of diagrams, it is often useful to observe which properties of the super-standard labeling are preserved by Kohnert moves. This provides a means for detecting whether $D \in \mathcal{P}(D_0)$ for some initial diagram $D_0$ without needing to explicitly check candidate sequences of Kohnert moves. This motivates the following definition:

\begin{definition}\label{def:NElabeling} 
    Let $D$ be a diagram. A \textbf{northeast labeling} $\mathcal{L}: D \to \mathbb{N}$ is a labeling that satisfies the following properties:
    \begin{enumerate}[label=(\arabic*)]
        \item\label{enum:astrict} $\mathcal{L}$ is strict.
        \item\label{enum:bbound} Each label in row $i$ is at least $i$.
        \item\label{enum:cnortheast} If cell $x'$ is to the right of cell $x$ and $\mathcal{L}(x') < \mathcal{L}(x)$, then there is a cell $x''$ in the same column as $x'$ with $\mathcal{L}(x) = \mathcal{L}(x'')$. 
        \item\label{enum:dbelow} If cell $x'$ is to the right of cell $x$ and $\mathcal{L}(x') = \mathcal{L}(x)$, then $x'$ is weakly below $x$.
    \end{enumerate}
    The pair $T = (D, \mathcal{L})$ is a \textbf{northeast tableau}.
\end{definition}

Our northeast tableaux are a generalization of the lock tableaux of Assaf and Searles \cite[Definition 6.5]{AS22}. The following lemma is used in the proofs of some of our main results.

\begin{lemma}\label{lemma:above}
    Let $D$ be a diagram with a northeast labeling $\mathcal{L}$. Let $x_1 = (r_1, c_1)$ and $x_2 = (r_2, c_2)$ be cells in $D$ such that $c_1 \leq c_2$ and $\mathcal{L}(x_1) > \mathcal{L}(x_2)$. Then $r_1 > r_2$.
\end{lemma}

\begin{proof}
    If $x_2$ and $x_1$ are in the same column, the result follows from the strictness of $\mathcal{L}$. Now suppose that $x_2 = (r_2, c_2)$ is strictly to the right of $x_1 = (r_1, c_1)$. By Property \ref{enum:cnortheast} of Definition \ref{def:NElabeling}, there is a cell $x_2' = (r_2', c_2)$ in $D$ such that $\mathcal{L}(x_2') = \mathcal{L}(x_1)$. Then $r_2' > r_2$ by strictness, and $r_1 \geq r_2'$ by Property \ref{enum:dbelow} of Definition \ref{def:NElabeling}. Therefore, $r_1 > r_2$.
\end{proof}
\noindent In other words, in a northeast tableau, the row indices of the string of cells labeled with an $i$ are weakly descending from left to right. This also means that the labels on the cells in a given row are weakly increasing from left to right.

Northeast tableaux capture northeast diagrams in the following sense: 

\begin{lemma}
\label{lemma:NEDiagramEquivSSLabelingNE}
    A diagram $D$ is northeast if and only if the super-standard labeling $\mathcal{L}$ on $D$ is a northeast labeling.
\end{lemma}

\begin{proof}
    First note that both northeast labelings and super-standard labelings always satisfy Properties \ref{enum:astrict}, \ref{enum:bbound}, and \ref{enum:dbelow} of Definition \ref{def:NElabeling}. 
    Therefore, it is sufficient to compare Property \ref{enum:cnortheast}. 
    It is straightforward to check that Property \ref{enum:cnortheast} being satisfied on a super-standard labeling is equivalent to the underlying diagram being northeast.
\end{proof}

\noindent Lemma \ref{lemma:NEDiagramEquivSSLabelingNE} motivates defining the \textbf{initial northeast tableau} of a northeast tableau $T$:

\begin{definition}
Given a northeast tableau $T = (D, \mathcal{L})$, the \textbf{initial northeast diagram} of $T$ is
\[D_0 \coloneqq \{(\mathcal{L}(x), c) \mid x = (r, c) \in T\}.\]
The \textbf{initial northeast tableau} $T_0 = (D_0, \mathcal{L}_0)$ of $T$ is the initial northeast diagram paired with the \textbf{initial northeast labeling} $\mathcal{L}_0: (r, c) \mapsto r$.
\end{definition}

\begin{proposition}
    The initial northeast tableau $T_0 = (D_0, \mathcal{L}_0)$ of a northeast tableau $T$ is the unique northeast tableau satisfying the following properties:
    \begin{enumerate}
        \item $\mathcal{L}_0$ is the super-standard labeling,
        \item $T_0$ is column-equivalent to $T$, and
        \item $D_0$ is a northeast diagram.
    \end{enumerate}
\end{proposition}

\begin{proof}
    Certainly there is at most one tableau endowed with the super-standard labeling that is column-equivalent to $T$ - its existence comes from the fact that $\mathcal{L}$ is strict, and thus injective when restricted to an individual column. Any super-standard labeling satisfies Properties \ref{enum:astrict}, \ref{enum:bbound}, and \ref{enum:dbelow} of Definition \ref{def:NElabeling}, and since $\mathcal{L}_0$ is column-equivalent to $\mathcal{L}$, it must also satisfy Property \ref{enum:cnortheast}. Thus, $T_0$ is a northeast tableau. Then $D_0$ is northeast by Lemma \ref{lemma:NEDiagramEquivSSLabelingNE}.
\end{proof}

Next, we introduce the concept of \textbf{total displacement}, which quantifies how far a labeling is from being super-standard. While defined in generality for Kohnert tableaux, the total displacement is well-suited for northeast tableaux because it provides a rank function to complement our result on classifying ranked northeast posets. 

\begin{definition}
    Let $T=(D, \mathcal{L})$ be a tableau. The \textbf{displacement} is a function $\delta_\mathcal{L}: D \to \mathbb{Z}_{\geq 0}$ such that for $x = (r,c) \in D$,  
    \[\delta_\mathcal{L}(x) = \mathcal{L}(x) - r.\]
    The \textbf{total displacement} of $T$ is
    \[\Delta_\mathcal{L}(T) = \sum_{x \in D} \delta_\mathcal{L}(x).\]
\end{definition}

See Figure \ref{fig:northeast} for an example.
If $T = (D, \mathcal{L})$ is a northeast tableau, we often assume that the labeling is standard with respect to the initial northeast diagram of $T$ and omit the subscript $\mathcal{L}$ from $\Delta$. 

\begin{figure}[ht] \centering
\begin{subfigure}{0.45\textwidth} \centering
    \begin{tikzpicture}[scale = 0.5]
    \draw (1,6)--(1,0)--(5,0);
    \draw (2,0) rectangle (3,1); \node at (2.5, 0.5) {$1$};%
    \draw (4,0) rectangle (5,1); \node at (4.5, 0.5) {$5$};%
    \draw (2,1) rectangle (3,2); \node at (2.5, 1.5) {$3$};%
    \draw (3,1) rectangle (4,2); \node at (3.5, 1.5) {$5$};%
    \draw (2,2) rectangle (3,3); \node at (2.5, 2.5) {$4$};%
    \draw (4,2) rectangle (5,3); \node at (4.5, 2.5) {$6$};%
    \draw (1,3) rectangle (2,4); \node at (1.5, 3.5) {$4$};%
    \draw (3,5) rectangle (4,6); \node at (3.5, 5.5) {$6$};%
    \end{tikzpicture}
    \caption{Figure \ref{subfig:northeast}. A northeast tableau $T$ with $\Delta_\mathcal{L}(T) = 12$.} \label{subfig:northeast}
\end{subfigure}
\begin{subfigure}{0.5\textwidth} \centering
    \begin{tikzpicture}[scale = 0.5]
    \draw (1,6)--(1,0)--(5,0);
    \draw (2,0) rectangle (3,1); \node at (2.5, 0.5) {$1$};%
    \draw (2,2) rectangle (3,3); \node at (2.5, 2.5) {$3$};%
    \draw (1,3) rectangle (2,4); \node at (1.5, 3.5) {$4$};%
    \draw (2,3) rectangle (3,4); \node at (2.5, 3.5) {$4$};%
    \draw (3,4) rectangle (4,5); \node at (3.5, 4.5) {$5$};%
    \draw (4,4) rectangle (5,5); \node at (4.5, 4.5) {$5$};%
    \draw (3,5) rectangle (4,6); \node at (3.5, 5.5) {$6$};%
    \draw (4,5) rectangle (5,6); \node at (4.5, 5.5) {$6$};%
    \end{tikzpicture}
    \caption{Figure \ref{subfig:northeastinitial}. The initial tableau $T_0$ of $T$ from Figure \ref{subfig:northeast}.} \label{subfig:northeastinitial}
\end{subfigure}
\caption{Two northeast tableaux.}\label{fig:northeast}
\end{figure}

\begin{lemma}\label{lemma:northeastPreserved}
   Let $T = (D, \mathcal{L})$ be a northeast tableau, and let $T_0 = (D_0, \mathcal{L}_0)$ be the initial northeast tableau of $T$. If $D'$ is a diagram obtained by performing a single Kohnert move on $D$, then the standard labeling $\mathcal{L}'$ of $D'$ with respect to $D_0$ is a northeast labeling.
\end{lemma}

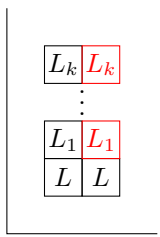
\begin{figure}[ht] \centering
    \begin{tikzpicture}[scale = 0.5]
    \draw (1,6)--(1,0)--(5,0);
    \draw (2,1) rectangle (3,2); \node at (2.5, 1.5) {$L$};%
    \draw (4,1) rectangle (5,2); \node at (4.5, 1.5) {$L$};%
    \draw (2,2) rectangle (3,3); \node at (2.5, 2.5) {$L_1$};%
    \draw[color=red] (4,2) rectangle (5,3); \node[color=red] at (4.5, 2.5) {$L_1$};%
    \node at (2.5, 3.7) {$\vdots$};
    \node at (3.5, 4.5) {$\hdots$};
    \draw (2,4) rectangle (3,5); \node at (2.5, 4.5) {$L_k$};%
    \draw[color=red] (4,4) rectangle (5,5); \node[color=red] at (4.5, 4.5) {$L_k$};%
    \node at (2.5, -0.6) {$c$};
    \node at (4.5, -0.5) {$c'$};
    \node at (0.5, 4.5) {$r$};
    \node at (0.5, 1.5) {$r'$};
    \end{tikzpicture}
    \caption{The diagram $D$ described in Step \ref{enum:prf4} of the proof of Lemma \ref{lemma:northeastPreserved}.} \label{fig:cond4}
\end{figure}

\begin{proof}
    By Lemma~\ref{lem:stdLabelingsConstruction}, the standard labeling of $D'$ may be constructed from the standard labeling for $D$. Observe that a cell $x'$ in $D'$ with label $n$ is weakly below and at most one position below the corresponding cell $x$ with the same label in the same column in $D$. We show that $\mathcal{L}'$ satisfies the four properties of Definition \ref{def:NElabeling} and hence is northeast.
    \begin{enumerate}[label=(\arabic*)]
        \item $\mathcal{L}'$ is strict by construction.
        \item Each label in $\mathcal{L}'$ in row $i$ is at least $i$ if and only if $\delta_{\mathcal{L}'}(x') \geq 0$ for all $x' \in D'$. If $x' \in D'$, then there exists an $x \in D$ in the same column with the same label. Then, $\delta_{\mathcal{L}'}(x') \geq \delta_\mathcal{L}(x)$ by our observation. However, $\delta_\mathcal{L}(x) \geq 0$ since $\mathcal{L}$ is northeast. Hence, $\delta_{\mathcal{L}'}(x') \geq 0$.
        
        \item Property \ref{enum:cnortheast} holds for $\mathcal{L'}$ because $\mathcal{L}$ and $\mathcal{L'}$ are column-equivalent, and $\mathcal{L}$ is northeast.
        
        \item \label{enum:prf4} Let $x = (r, c)$ be the cell to which the Kohnert move on $D$ was applied. For the sake of contradiction, suppose that $\mathcal{L}'$ does not satisfy Property \ref{enum:dbelow}. Then, there is a pair of cells $y_1', y_2' \in D'$ such that $\mathcal{L}'(y_1') = \mathcal{L}'(y_2')$, but $y_1'$ is to the left of and strictly below $y_2'$. Because the only difference between $D$ and $D'$ is in column $c$, it follows that column $c$ contains $y_1'$ in a position weakly below $x$. Moreover, the cell $y_1$ in column $c$ of $D$ with the label $\mathcal{L}'(y_1')$ is exactly one position above $y_1'$. Then because $\mathcal{L}$ is northeast, $y_1$ is in the same row as $y_2'$. Setting $y_1' = (r' - 1, c)$ and $y_2' = (r', c')$ for some $r' \leq r$ and $c' > c$, we have that $y_1 = (r', c)$ and $y_2 = (r', c')$ are the corresponding cells in $D$ with \[\mathcal{L}(y_1) = \mathcal{L}'(y_1') = \mathcal{L}'(y_2') = \mathcal{L}(y_2).\] Denote this label by $L := \mathcal{L}(y_1)$. Because $y_1'$ is below $x$, $y_1$ is weakly below $x$. Then by strictness, $L \leq \mathcal{L}(x)$. Because $y_1$ was affected by the Kohnert move, there are no empty positions between $x$ and $y_1$ in $D$. In particular, there are $r - r' + 1$ cells between $x$ and $y_1$ in $D$, including $x$ and $y_1$. 
        
        If $r=r'$, then $x$ is not rightmost in its row and thus can't be moved, which is a contradiction. Thus, we can suppose $x\neq y_1$ and let $z_1$ denote the cell directly above $y_1$ in $D$. Since $\mathcal{L}$ is a northeast labeling, there is a cell $z_2$ in column $c'$ that is weakly below $z_1$ but strictly above $y_2$ with the same label as $z_1$. Then, $z_2$ is in the same row as $z_1$. Similar logic shows that for each cell in column $c$ between $x$ and $y_1$ in $D$, including $x$ and $y_1$, there is a corresponding cell with the same label in the same row of column $c'$ in $D$. This contradicts $x$ being the cell to which the Kohnert move was applied, because $x$ is not the rightmost cell in its row.    

        Figure \ref{fig:cond4} depicts what the diagram $D$ must have looked like at the conclusion of Step \ref{enum:prf4}. The red cells are the ones whose existence is inferred by the northeast property of the labeling. The cell $x$ is the top left cell, labeled $L_k$. It clearly cannot be the subject of a Kohnert move, as it is not the rightmost cell in its row.
    \end{enumerate} 
    Since $\mathcal{L}'$ satisfies the four properties of Definition~\ref{def:NElabeling}, it is a northeast labeling.
    \end{proof}

\begin{remark} \label{rmk:northeast_labeling}
    Using Lemma \ref{lemma:northeastPreserved} and induction, we see that if $D \leq D_0$ and $D_0$ is a northeast diagram, then the standard labeling of $D$ with respect to $D_0$ is a northeast labeling.
\end{remark}

\begin{lemma}\label{lemma:incrementDisplacement}
    Let $T = (D, \mathcal{L})$ be a northeast tableau with nonzero total displacement. Then, there exists a northeast tableau $T' = (D', \mathcal{L}')$ that is column-equivalent to $T$ and such that $D<D'$ and $\Delta_{\mathcal{L}'}(T') = \Delta_\mathcal{L}(T) - 1$.
\end{lemma}

\begin{proof}
    Because $\Delta_\mathcal{L}(T) > 0$, $D$ has at least one cell with nonzero displacement. Let $k$ be the maximum label among all cells with non-zero displacement and $x=(r,c)$ be the leftmost among such cells.
    We claim that $D$ has no cells $x' = (r + 1, c')$ with $c' \geq c$. Suppose for the sake of contradiction that such a cell $x'$ does exist. We examine this possibility in cases:

\begin{itemize}[leftmargin=0.65in]
\item[\textbf{Case 1:}] Suppose  $\mathcal{L}(x') > \mathcal{L}(x)$. Then $\delta_\mathcal{L}(x') = \mathcal{L}(x') - (r+1) \geq \mathcal{L}(x) - r = \delta_\mathcal{L}(x) > 0$, contradicting the maximality of the label of $k$ among all cells with nonzero displacement. 

\item[\textbf{Case 2:}] Suppose $\mathcal{L}(x') = \mathcal{L}(x)$. Since $\mathcal{L}$ is strict, we must have $c'>c$.  This contradicts Property \ref{enum:dbelow} of northeast tableaux.

\item[\textbf{Case 3:}] Suppose $\mathcal{L}(x') < \mathcal{L}(x)$. Again, since $\mathcal{L}$ is strict, we have $c'>c$. By Property \ref{enum:cnortheast} of northeast tableaux, there is a cell $x''$ in the same column as $x'$ with $\mathcal{L}(x) = \mathcal{L}(x'')$. Since $\mathcal{L}$ is strict, $x''$ is above $x'$, and so the pair $x,x''$ contradicts Property~\ref{enum:dbelow} of northeast tableaux.
\end{itemize} 

    Then $D$ has no cells $x' = (r + 1, c')$ with $c' \geq c$. Now say
    \[D' = (D \setminus \{(r, c')\}) \cup \{(r+1, c')\}\]
    and say $\tilde{x} \coloneqq (r+1, c') \in D'$. By our first observation, $\tilde{x}$ has no cells to its right in $D'$, and the position below it is empty. Therefore, $D < D'$ via a Kohnert move applied to $\tilde{x}$. Let $\mathcal{L}'$ be the unique strict labeling of $D'$ that is column-equivalent to $\mathcal{L}$. Note that $\mathcal{L}$ agrees with $\mathcal{L'}$ on the intersection of $D$ and $D'$, and the cells in each diagram not in that intersection ($x$ and $\tilde{x}$, respectively) have the same label. Thus, $\Delta_\mathcal{L}(T') = \Delta_\mathcal{L}(T) - 1$. We now check that $\mathcal{L}'$ satisfies the four properties of Definition \ref{def:NElabeling} and hence is northeast.
    \begin{enumerate}[label=(\arabic*)]
        \item $\mathcal{L}'$ is strict by construction.
        \item  Each label in $\mathcal{L}'$ in row $i$ is at least $i$ if and only if $\delta_{\mathcal{L}'}(x') \geq 0$ for all $x' \in D'$. If $x' \neq \tilde{x}$, then $\delta_{\mathcal{L}'}(x') = \delta_{\mathcal{L}}(x') > 0$ since $\mathcal{L}$ is northeast, while $\delta_{\mathcal{L}'}(\tilde{x}) = \delta_{\mathcal{L}}(x) - 1 \geq 0$  by the choice of $x$ and the definition of $\tilde{x}$.
        \item Property \ref{enum:cnortheast} holds for $\mathcal{L'}$ because $\mathcal{L}$ and $\mathcal{L'}$ are column-equivalent and $\mathcal{L}$ is northeast.
        \item  For the sake of contradiction, suppose that $\mathcal{L}'$ does not satisfy Property \ref{enum:dbelow}. Since $\mathcal{L}$ and $\mathcal{L'}$ agree on the intersection of $D$ and $D'$, and  $\mathcal{L}$ is northeast, any violation of Property \ref{enum:dbelow} must occur for a pair of cells with the rightmost cell equal to $\tilde{x}$. Let $x'$ be the leftmost cell in the pair. Since $x'$ and $\tilde{x}$ violate Property \ref{enum:dbelow}, $x'$ is strictly below $\tilde{x}$ in $D'$, and hence weakly below $x$ in $D$. Because $\mathcal{L}$ satisfies Property \ref{enum:dbelow}, $x'$ cannot be strictly below $x$. Then $x'$ and $x$ are in the same row, so $\delta_{\mathcal{L}}(x') = \delta_{\mathcal{L}}(x)$. This contradicts $x$ being the leftmost cell in its row with nonzero displacement. We conclude that Property \ref{enum:dbelow} holds for $\mathcal{L}'$. \qedhere
    \end{enumerate}
\end{proof}

\begin{theorem}\label{thm:labeling}
    Let $D_0$ be a northeast diagram with super-standard labeling $\mathcal{L}_0$. Then $D \in \mathcal{P}(D_0)$ if and only if $D$ admits a northeast labeling $\mathcal{L}$ that is column-equivalent to $\mathcal{L}_0$. 
\end{theorem}

\begin{proof}
    (Proof of $\Rightarrow$)  Suppose $D \in \mathcal{P}(D_0)$. By Remark~\ref{rmk:northeast_labeling}, we can inductively apply Lemma~\ref{lemma:northeastPreserved} to obtain a northeast labeling $\mathcal{L}$ of $D$ that is a standard labeling with respect to $D_0$. Thus, $\mathcal{L}$ is column-equivalent to $\mathcal{L}_0$. 
    
    (Proof of $\Leftarrow$) Let $T = (D, \mathcal{L})$ be a northeast tableau whose labeling $\mathcal{L}$ is column-equivalent to $\mathcal{L}_0$. We induct on $\Delta_\mathcal{L}(T)$. If $\Delta_\mathcal{L}(T) = 0$ then $\mathcal{L}$ is the super-standard labeling, so $D = D_0$. Now suppose that $\Delta_\mathcal{L}(T) > 0$ and for all $T' = (D', \mathcal{L}')$ such that $T$ and $T'$ are column-equivalent and $\Delta_{\mathcal{L}'}(D') = \Delta_\mathcal{L}(D) - 1$, we have $D' \in \mathcal{P}(D_0)$. By Lemma \ref{lemma:incrementDisplacement}, there is at least one such $T'$ with $D < D'$. Then, by transitivity $D \in \mathcal{P}(D_0)$.
\end{proof}

\subsection{Elementary chains}
We begin this section with a proof of Theorem \ref{thm:NEKohnertPosetEle}. We then proceed to show that for any saturated chain in a Kohnert poset that is defined by the application of solely elementary moves, the moves can always be reordered so that cells labeled by smaller numbers are always moved first.

\begin{lemma}\label{lem:NEavoidsJumps}
    Let $D_0$ be a northeast diagram and let $D \in \mathcal{P}(D_0)$. Then $D$ can be obtained from $D_0$ be a (possibly empty) sequence of elementary Kohnert moves.
\end{lemma}

\begin{proof}
    Let $\mathcal{L}_0$ be the super-standard labeling of $D_0$. Then $T_0 = (D_0,\mathcal{L}_0)$ is a northeast tableau by Lemma \ref{lemma:NEDiagramEquivSSLabelingNE}. By Theorem \ref{thm:labeling}, there exists a labeling $\mathcal{L}$ such that $T = (D,\mathcal{L})$ is a northeast tableau that is column-equivalent to $\mathcal{L}_0$.
    
    We induct on $\Delta_{\mathcal{L}}(T)$. If $\Delta_{\mathcal{L}}(T) = 0$, then $D = D_0$ and the desired sequence of elementary moves is the empty sequence. Now, fix a diagram $D$ with $\Delta_{\mathcal{L}}(T) = d$ and $d > 0$. 
    By Lemma \ref{lemma:incrementDisplacement}, there is a diagram $D'$ with column-equivalent tableau $T'=(D',\mathcal{L}')$ such that $D < D'$ and $\Delta_{\mathcal{L'}}(T') = \Delta_{\mathcal{L}}(T) - 1$. Theorem \ref{thm:labeling} implies that $D' \in \mathcal{P}(D_0)$. Hence, by induction, there is a sequence of elementary Kohnert moves connecting $D_0$ to $D'$ in $\mathcal{P}(D_0)$. The Kohnert move that connects $D$ and $D'$ is elementary by construction. Thus, there is a sequence of elementary moves to get from $D_0$ to $D$ in $\mathcal{P}(D_0)$. 
\end{proof}

Our first main theorem follows from Lemma~\ref{lem:NEavoidsJumps}. 

\begin{repeattheorem}[Theorem 1.2]
 Let $D$ be a northeast diagram. Then, $\mathcal{P}(D)$ is a refinement of $\mathcal{P}^{ele}(D)$; in particular, $\mathcal{P}(D)$ and $\mathcal{P}^{ele}(D)$ are equal as sets.
\end{repeattheorem}

\begin{proof}[Proof of Theorem~\ref{thm:NEKohnertPosetEle}]
Let $D_1, D_2 \in \mathcal{P}(D)$ such that $D_1 \leq_{\mathcal{P}(D)} D_2$. Then it is clear that $D_1 \in \mathcal{P}(D_2)$. By Lemma \ref{lem:NEavoidsJumps}, $D_1$ can be obtained from $D_2$ by a sequence of elementary Kohnert moves. Hence $D_1 \leq_{\mathcal{P}^{ele}(D)}(D_2)$. This proves that $\mathcal{P}(D)$ is a refinement of $\mathcal{P}^{ele}(D)$. 

By definition, $\mathcal{P}^{ele}(D) \subseteq \mathcal{P}(D)$ as sets. Since $\mathcal{P}(D)$ is a refinement of $\mathcal{P}^{ele}(D)$, we have $\mathcal{P}^{ele}(D) = \mathcal{P}(D)$ as sets.
\end{proof}

The final definitions and results of this section are useful for proving our theorem on boundedness for northeast diagram posets.
\begin{definition}\label{def:elementary_chain}
Let $D_0$ be a diagram. A saturated chain in $\mathcal{P}(D_0)$, denoted by $D_1 \gtrdot D_2 \gtrdot \cdots \gtrdot D_M$, is an \textbf{elementary chain} if the Kohnert move connecting $D_a$ to $D_{a+1}$ is an elementary Kohnert move for all $1 \leq a < M$. 
\end{definition}

\begin{definition}\label{def:move_inversion}
Let $T_0 = (D_0, \mathcal{L}_0)$ be a tableau, and let $\mathbf{D} = \{D_1 \gtrdot D_2 \gtrdot \cdots \gtrdot D_M\}$ be an elementary chain. For any $D_a\gtrdot D_{a+1}$ in $\mathbf{D}$, define $\movelabel(D_a,D_{a+1})$ to be $i$ when $D_a$ and $D_{a+1}$ are connected by an elementary Kohnert move affecting a cell with label $i$. Define
\[
\mathcal{I}(\mathbf{D}) = \{ (a,b) \mid 1 \leq a < b < M\text{ and }\movelabel(D_a,D_{a+1}) > \movelabel(D_b,D_{b+1}) \}
\]
and $i(\mathbf{D})=|\mathcal{I}(\mathbf{D})|$. The pair $(a, b)$ is said to be a \textbf{move inversion}.
\end{definition}

The use of the word ``inversion'' is motivated by the fact that a move inversion $(a, b)$ represents a pair of moves in the chain where a larger-labeled cell has moved before a smaller-labeled cell. The following lemma shows that any elementary chain with maximal element $D_0$ may be replaced by an elementary chain with identical maximal and minimal elements that contains no move inversions. We highlight that this result holds for all diagrams. 

\begin{lemma}
\label{lemma:elementaryChainShuffle}
Let $T_0 = (D_0,\mathcal{L}_0)$ be a tableau with the super-standard labeling, and $\mathbf{D} := D_1 \gtrdot \cdots \gtrdot D_M$ be an elementary chain in $\mathcal{P}(D_0)$. 
Then, there exists an elementary chain $\mathbf{F} := F_1 \gtrdot \cdots \gtrdot F_M$ such that $F_1 = D_1$, $F_M=D_M$, and $i(\mathbf{F})=0$.
\end{lemma}

\begin{proof}
If $i(\mathbf{D})=0$, then set $\mathbf{F}=\mathbf{D}$. Otherwise, if $i(\mathbf{D}) > 0$, we will show that a new chain $\mathbf{E} = E_1 \gtrdot \cdots \gtrdot E_M$ may be constructed such that $D_1 = E_1$, $D_M = E_M$, and $i(\mathbf{E}) < i(\mathbf{D})$. This proves our desired result since, as the number of move inversions decreases at each step of the process, the iteration of this process must terminate.

Let $i(\mathbf{D}) > 0$, and take $(a,b) \in \mathcal{I}(\mathbf{D})$ such that $(a,b)$ is lexicographically minimal. The minimality of $(a,b)$ implies that $\movelabel(D_a,D_{a+1})\leq \movelabel(D_k,D_{k+1})\text{ for }a < k < b$ which in turn implies
\begin{equation}
\label{eq:shuffleq1}
\movelabel(D_k,D_{k+1}) > \movelabel(D_b,D_{b+1})\text{ for }a < k < b.
\end{equation}

For each $a \leq k \leq b$, let $x_k=(r_k, c_k)$ be the cell that is moved from $D_k$ to $D_{k+1}$. Define $E_{a+1}$ to be the diagram obtained from $D_a$ by the elementary Kohnert move of the cell $x_b = (r_b, c_b)$. We first verify that this move is valid. The cell $x_b$ exists in $D_a$ by \eqref{eq:shuffleq1}. Additionally, suppose, for contradiction, that $x_b$ is not the rightmost cell in its row in $D_a$. Since $x_b$ is the rightmost cell in its row in $D_b$, by \eqref{eq:shuffleq1}, we must have that the rightmost cell $x'$ in row $r_b$ in $D_a$ has label $\ell'$ for some $\ell' > \movelabel(D_b,D_{b+1})$. Since $\mathcal{L}_0$ is super-standard, the cell $x'$ labeled $\ell'$ could not have started in row $r_b$. So there exists some $1 \leq q < a$ such that $\movelabel(D_q,D_{q+1})=\ell'$. This contradicts the minimality of $(a,b)$. Likewise, because $(r_b - 1, c_b) \not\in D_a$ but $(r_b, c_b) \in D_a$, any move that would place a cell in some $D_k$ in position $(r_b -1 , c_b)$ without affecting cell $x_b$ could not be elementary. Thus, the move from $D_a$ to $E_{a+1}$ is a valid elementary Kohnert move.

For $a < k \leq b$, iteratively define $E_{k+1}$ to be the diagram reached from $E_{k}$ by the elementary Kohnert move of the cell $(r_{k-1}, c_{k-1})$. We now verify that each of these moves is valid. The cell $(r_{k-1}, c_{k-1})$ exists in $E_{k}$ because the diagram $E_{k}$ differs from $D_{k-1}$ only by the elementary Kohnert move of the cell $x_b=(r_b, c_b)$, which was moved first. Additionally, suppose, for contradiction, that $(r_{k-1}, c_{k-1})$ is not the rightmost cell in its row in $E_{k}$. Since $(r_{k-1}, c_{k-1})$ is the rightmost cell in its row in $D_{k-1}$, it must be the case that the rightmost cell in $E_{k}$ is the $\movelabel(D_b,D_{b+1})$-cell. Thus, row $r_{k-1}$ in $E_k$ contains both a $\movelabel(D_b,D_{b+1})$-cell and a $\movelabel(D_{k-1},D_{k})$-cell. Then \eqref{eq:shuffleq1} implies that $\movelabel(D_{k-1},D_{k}) > \movelabel(D_b,D_{b+1})$. Hence, there exists some $1 \leq q < k - 1$ such that both $\movelabel(D_q,D_{q+1})=\movelabel(D_{k-1},D_{k})$ and the Kohnert move connecting $D_q$ and $D_{q+1}$ moves the $\movelabel(D_{k-1},D_{k})$-cell into row $r_{k-1}$. If $1 \leq q < a$, then this contradicts the minimality of $(a,b)$. Otherwise, if $a \leq q < k-1$, then $(r_q, c_q)$ is not the rightmost cell in its row in $D_q$; the rightmost cell is the $\movelabel(D_b,D_{b+1})$-cell $(r_b, c_b)$ (see Figure \ref{fig: pf helper}). This is a contradiction of the definition of $\mathbf{D}$. In either case, we have a contradiction, so we conclude that $(r_{k-1}, c_{k-1})$ is the rightmost cell in its row in $E_{k}$. Therefore, for each $a < k \leq b$, the move from $E_k$ to $E_{k+1}$ is a valid elementary Kohnert move.

For $1 \leq k \leq a$ or $b+1 < k \leq M$, we define $E_k = D_k$. Note that $E_a = D_a$ and $E_{b+1}=D_{b+1}$. Then $\mathbf{E} = E_1 \gtrdot E_2 \gtrdot \cdots \gtrdot E_M$ is a new chain such that $D_1 = E_1$, $D_M = E_M$, and $i(\mathbf{E}) = i(\mathbf{D}) - (b-a) < i(\mathbf{D})$.
\end{proof}

\begin{center}
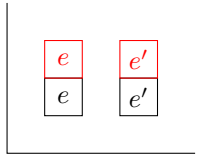

    \begin{tikzpicture}[scale=0.5]
    \celx(1,1)[e] \celx(3,1)[e'] 
    \rcelx(1,2)[e]    \rcelx(3,2)[e']
    \draw (0,4)--(0,0)--(5,0);
    \node at (1.5, -0.5) {$c_{k-1}$};
    \node at (-0.8, 1.5) {$r_{k-1}$};
    \end{tikzpicture}
    \captionof{figure}{If $(r_{k-1}, c_{k-1})$ with label $e=\movelabel(D_{k-1},D_k)$ were not rightmost in $E_k$ (black boxes), then it would need to be blocked by a box with label $e'=\movelabel(D_b,D_{b+1})$. But then this cell would not be rightmost when it is at $(r_{k-1}+1,c_{k-1})$ in $D_q$ (red boxes) in the case of $a\leq q \leq k-1$.}
    \label{fig: pf helper}
\end{center}

Combining Lemmas \ref{lem:NEavoidsJumps} and \ref{lemma:elementaryChainShuffle}, we have the following proposition which plays a key role in our proof of Theorem \ref{thm:mainNEbounded}.

\begin{proposition}\label{prop:NE elem chain no inversions}
If $D_0$ is northeast and $D \in \mathcal{P}(D_0)$, then there exists an elementary chain with maximal element $D_0$ and minimal element $D$ that does not contain any move inversions.
\end{proposition}

\section{Monomial multiplicity-free polynomials}\label{sec:Mfree}

We now prove our classification of monomial multiplicity-free Kohnert polynomials of northeast diagrams. We start by restating the theorem for convenience. We refer the reader to Figure \ref{fig:mmf_sub1} for a visualization of conditions (a) - (d) and to Figure~\ref{fig:mmf_sub2} for a northeast diagram whose Kohnert poset is not multiplicity-free.

\begin{figure}[h]
\centering
\begin{subfigure}[t]{.55\textwidth}
  \centering
  \begin{tikzpicture}[scale=0.5]
    \celx(1,2)[x_1]
    \celx(3,4)[x_2]
    \draw[pattern=dots] (1,0) rectangle (2,2);
    \draw[pattern=vertical lines] (2,0) rectangle (5,3);
    \draw (0,5)--(0,0)--(5,0)--(5,5);
    \node at (1.5,-0.5) {$c_1$};
    \node at (3.5,-0.5) {$c_2$};
    \node at (-0.5,2.5) {$r_1$};
    \node at (-0.5,4.5) {$r_2$};
    \end{tikzpicture}
  \caption{Figure~\ref*{fig:mmf_sub1}. The arrangement of cells $x_1$ and $x_2$ satisfying (a) and (b). Condition (c) states that there is an empty position in the dotted region. Condition (d) states that there are at least two empty positions in each column of the striped region. }
  \label{fig:mmf_sub1}
\end{subfigure}\qquad
\begin{subfigure}[t]{.4\textwidth}
  \centering
  \begin{tikzpicture}[scale = 0.5]
    \draw (1,4)--(1,0)--(3,0);
    \cel(1,1)
    \cel(2,2)
    \end{tikzpicture}
  \caption{Figure~\ref*{fig:mmf_sub2}. A minimal example of a northeast diagram whose Kohnert poset does not yield a monomial multiplicity-free polynomial.}
  \label{fig:mmf_sub2}
\end{subfigure}
\end{figure}

\begin{repeattheorem}[Theorem~\ref{thm:mainMFree}]
    If $D_0$ is a northeast diagram, then $\mathfrak{K}_{D_0}$ is monomial multiplicity-free if and only if $D_0$ does not contain $x_1 = (r_1,c_1)$ and $x_2 = (r_2,c_2)$ such that:
    \begin{itemize}
        \item[(a)] $r_1 < r_2$,
        \item[(b)] $c_1 < c_2$,
        \item[(c)] there is at least one empty position $(r, c_1)$ where $r < r_1$, and
        \item[(d)] for each $c>c_1$, there are at least two empty positions $(r,c)$ where $r \leq r_1$. 
    \end{itemize}
\end{repeattheorem}

\begin{proof}[Proof of Theorem~\ref{thm:mainMFree}]
    (Proof of $\Rightarrow$) We prove the contrapositive. Suppose $D_0$ has a pair of cells $x_1 = (r_1, c_1)$ and $ x_2 = (r_2, c_2)$ with the described properties. If $D_0$ has multiple such pairs, choose the pair such that the word $(c_2, -r_2, c_1, r_1)$ is maximal in the lexicographic order on $\Z^4$. By this choice of $x_1$ and $x_2$, the following are true:
    \begin{itemize}
        \item The only cell $(r, c)$ such that $r_1 \leq r < r_2$ and $c_1 \leq c < c_2$ is $x_1$.
        \item The only cell $(r, c)$ such that $r = r_2$ and $c \geq c_2$ is $x_2$.
        \item If present, $(r_1,c_2)$ must have at least two empty positions below it.
        \item The only possible cell $(r, c)$ such that $r_1 \leq r < r_2$ and $c \geq c_2$ is $(r_1, c_2)$. 
        \item Column $c_2$ must be the rightmost nonempty column in $D_0$.           
    \end{itemize}

    For the last two bullet points, the violation of either property, combined with the fact that $D_0$ is northeast, would imply the presence of a cell that contradicts our choice of $x_2$.

    We perform the following sequence of operations on $D_0$ to obtain a diagram $D_1$:
    \begin{enumerate}
        \item If $(r_1, c_2) \in D_0$ or $(r_1 - 1, c_2) \in D_0$, we perform Kohnert moves on $D_0$ to obtain a diagram where $(r_1, c_2)$ and $(r_1, c_2 - 1)$ are empty. 
        \item If $(r_1 - 1, c) \in D_0$ for any $c_1 < c \leq c_2$, we perform Kohnert moves to obtain a diagram where $(r_1 - 1, c)$ is empty. 
        \item By condition $(c)$, if there is a cell directly below $x_1$, we can move it as well, so that the position $(r_1-1,c_1)$ is empty.
    \end{enumerate}
    
    This creates the diagram $D_1$, where $x_1$ and $x_2$ are unmoved and the following are satisfied:
    \begin{itemize}
        \item The only cell $(r, c)$ such that $r_1 - 1 \leq r < r_2$ and $c \geq c_1$ is $x_1$.
        \item The only cell $(r, c)$ such that $r = r_2$ and $c \geq c_2$ is $x_2$.
    \end{itemize}

    In other words, $x_1$ is rightmost in $r_1$, and $x_2$ will remain rightmost in its row if it moves down to any row of index $\geq r_1-1$. 
    Then, the diagram $D_2$ obtained from $D_1$ by performing a single Kohnert move on $x_1$ and $r_2 - r_1$ Kohnert moves on $x_2$ has the same row weight as the diagram $D_2'$ obtained from $D_1$ by performing $r_2 - r_1 + 1$ Kohnert moves on $x_2$. 
    Thus, $\mathcal{P}(D_0)$ is not monomial multiplicity-free.
    
    \hfill
    
    (Proof of $\Leftarrow$) Suppose that $\mathfrak{K}_{D_0}$ is not monomial multiplicity-free. In other words, there exist distinct $D, D' \in \mathcal{P}(D_0)$ such that $\rwt(D) = \rwt(D')$. Let $\mathcal{L}_0, \mathcal{L}$, and $\mathcal{L}'$ denote the standard labelings on $D_0, D$, and $D'$, respectively.
    
    Let row $r$ denote the highest row where $D$ and $D'$ differ. Choose $c_1$ and $c_2$ to be the two maximum column indices such that $y_1 = (r, c_1)$ is a cell in $D \setminus D'$ and $y_2' = (r, c_2)$ is a cell in $D' \setminus D$. Without loss of generality, assume that $c_1 < c_2$. Define $L_1 = \mathcal{L}(y_1)$ and $L_2 = \mathcal{L}'(y_2')$. 
    
    Since $D$ and $D'$ are identical above row $r$, the column-equivalence and strictness of $\mathcal{L}$ and $\mathcal{L}'$ imply that all of the cells above row $r$ in $D$ have the same labels as the corresponding cells in $D'$. 
    Thus, the cells $y_1' = (r_1, c_1) \in D'$     
    and $y_2 = (r_2, c_2) \in D$ such that $\mathcal{L}'(y_1') = L_1$ and $\mathcal{L}(y_2) = L_2$ must have $r_1, r_2 < r$. See Figure \ref{fig:MMF helper}.

\begin{figure}[ht]
  \begin{center}
    \begin{tikzpicture}[scale=0.5]
    \celx(1,3)[y_1] \celx(3,2)[y_2] 
    \draw (0,5)--(0,0)--(5,0);
    \node at (2.5,-1.5) {$D$};
    \node at (1.5, -0.5) {$c_1$};
    \node at (3.5, -0.5) {$c_2$};
    \node at (-0.5, 3.5) {$r$};
    \node at (-0.5, 2.5) {$r_2$};
    \end{tikzpicture}
    \hspace{3em}
      \begin{tikzpicture}[scale=0.5]
    \celx(1,1)[y_1'] \celx(3,3)[y_2'] 
    \draw (0,5)--(0,0)--(5,0);
    \node at (2.5,-1.5) {$D'$};
    \node at (1.5, -0.5) {$c_1$};
    \node at (3.5, -0.5) {$c_2$};
    \node at (-0.5, 3.5) {$r$};
    \node at (-0.5, 1.5) {$r_1$};
    \end{tikzpicture}
    \hspace{3em}
    \begin{tikzpicture}[scale=0.5]
    \celx(1,3)[x_0] \celx(3,4)[x_2] \celx(2,3)[x_1]  
    \draw (0,5)--(0,0)--(5,0);
    \node at (2.5,-1.5) {$D_0$};
    \node at (1.5, -0.5) {$c_1$};
    \node at (2.5, -0.5) {$\tilde{c_1}$};
    \node at (3.5, -0.5) {$c_2$};
    \node at (-0.5, 3.5) {$L_1$};
    \node at (-0.5, 4.5) {$L_2$};
    \end{tikzpicture}
    \caption{The relative positions of cells $y_1, y_2$ in $D$ and $y_1', y_2'$ in $D'$, respectively, as well as the cells in $D_0$ with labels $L_1$ and $L_2$. 
    }\label{fig:MMF helper}
  \end{center}
\end{figure}
    Note that $y_1'$ is strictly below and to the left of $y_2'$. By Property \ref{enum:dbelow} of Definition \ref{def:NElabeling}, $L_1 \neq L_2$. 
    Moreover, $L_1 > L_2$ would contradict Lemma \ref{lemma:above}. 
    Therefore, $L_1 < L_2$.

    Define $x_2 = (L_2, c_2) \in D_0$. 
    Let $x_1 = (L_1, \tilde{c}_1)$ be the cell in $D_0$ such that $\tilde{c}_1 < c_2$ and $\tilde{c}_1$ is maximal. 
    Such a choice of cell exists because $x_0:=(L_1, c_1) \in D_0$ and $c_1 < c_2$. 
    Because $L_1 < L_2$ and $\tilde{c}_1 < c_2$, $x_1$ and $x_2$ satisfy Properties $(a)$ and $(b)$ of theorem statement \ref{thm:mainMFree}. 
    We now show that $x_1$ and $x_2$ satisfy Properties $(c)$ and $(d)$ as well.

    \begin{itemize}
        \item[(c)] Observe that $\mathcal{L}_0(x_1) = L_1$.
        In $D'$, the cell $y_1'$ has nonzero displacement since it occurs below the cell with label $L_1$ in the same column in $D$. Then, by Property \ref{enum:dbelow} of Definition \ref{def:NElabeling}, the cell labeled $L_1$ in column $\tilde{c}_1$ of $D'$ must also have nonzero displacement. Thus, there exists at least one empty space below $x_1$ in $D_0$.
        \item[(d)] If there are no cells to the right of and weakly below $x_1$ in $D_0$, then because $x_1$ is not in the first row, this property is satisfied. Now suppose that there is some cell $x_3 = (r_3, c_3)$ in $D_0$ that is to the right of and weakly below $x_1$. Since $D_0$ is northeast, there is a cell $x_4 = (L_1, c_3)$ in $D_0$ that may or may not be distinct from $x_3$. 
        By our choice of $\tilde{c}_1$, we have $c_3 \geq c_2$. Let $y_4 \in D$ denote the cell in column $c_3$ with label $L_1$. In other words, $y_4$ is where $x_4$ ended up after the sequence of Kohnert moves that produced $D$ from $D_0$. The cell $y_2$ is strictly below $y_1$ in $D$. And since $\mathcal{L}(y_2)=L_2 > L_1 = \mathcal{L}(y_4)$, Lemma \ref{lemma:above} implies that $y_4$ is strictly below $y_2$ in $D$. Then $\delta_{\mathcal{L}}(y_4) \geq \delta_{\mathcal{L}}(y_1) + 2 \geq 2$, so $y_4$ moved at least twice in the sequence of Kohnert moves used to obtain $D$ from $D_0$. Thus, in $D_0$, $x_4$ has at least two empty spaces below it. 
        Therefore, all columns $c >\tilde{c}_1$ in $D_0$ must have at least two empty positions in rows below or equal to $L_1$ as desired. \qedhere
    \end{itemize}

\end{proof}

\section{Ranked posets}\label{sec:ranked}

Before proving our criterion for rankedness of a Kohnert poset of a northeast diagram, we first note that the total displacement function $\Delta$ described in Section~\ref{sec:labelings} can be used to construct a rank function when applicable.

\begin{lemma}\label{lem:rankedCovering}
    Let $D_0$ be a northeast diagram with the super-standard labeling $\mathcal{L}_0$, and let $T_0=(D_0,\mathcal{L}_0)$. Assign all $D\in \mathcal{P}(D_0)$ the standard labeling. Then $\mathcal{P}(D_0)$ is ranked if and only if it has $-\Delta$ as a rank function. Equivalently, every covering relation in $\mathcal{P}(D_0)$ is given by an elementary Kohnert move.
\end{lemma}

\begin{proof}
    First, if $\mathcal{P}(D_0)$ has $-\Delta$ as a rank function, then by definition, $\mathcal{P}(D_0)$ is ranked. 
    
    Next, suppose that $\mathcal{P}(D_0)$ is ranked. By Lemma \ref{lemma:incrementDisplacement}, if $D \in \mathcal{P}(D_0)$ has nonzero displacement, then there is some $D' \in \mathcal{P}(D_0)$ with $D < D'$ and $\Delta(T) = \Delta(T') + 1$. Note that two diagrams with equal displacement are not comparable, since comparable diagrams differ by a sequence of Kohnert moves, each of which increase displacement. So, there is no $D''$ such that $D < D'' < D'$. Thus, $D$ is covered by $D'$.
    
    By induction, there is a maximal chain of length $\Delta(T)$ in $\mathcal{P}(D_0)$ that connects $D_0$ to $D$. Since $\mathcal{P}(D_0)$ is ranked, all maximal chains from $D_0$ to $D$ have length $\Delta(T)$. As applying a Kohnert move always increases displacement, each covering relation in each maximal chain changes displacement by 1. Thus, $\mathcal{P}(D_0)$ is ranked if and only if it has $-\Delta$ as a rank function.
\end{proof}

Next, we record a result of Colmenarejo, Hutchins, Mayers, and Phillips \cite{CHMP23} on general ranked Kohnert posets which we apply in the proof of Theorem~\ref{thm:mainNEranked}.
\begin{lemma}[{\cite[Theorem 3.4]{CHMP23}}]\label{lem:chmp}
    Let $D_0$ be a diagram. Suppose that there exists a diagram $D \in \mathcal{P}(D_0)$ and natural numbers $r, c_1, c_3 \in \mathbb{N}$ (with $1 \leq c_1 < c_3$) satisfying either conditions (a)(i) through (a)(iv) or conditions (b)(i) through (b)(iv):
    \begin{enumerate}[label=(\alph*)]
    \item
    \begin{itemize}
        \item[(i)] $(r +1, c_1), (r+2, c_1), (r+1, c_3) \in D$.
        \item[(ii)] $(r+2, c') \not\in D$ for all $c' > c_1$.
        \item[(iii)] $(r+1, c')\not\in D$ for all $c' > c_1$ not equal to $c_3$.
        \item[(iv)] $(r, c_1), (r, c_3) \not\in D$.
    \end{itemize}
    \item
        \begin{itemize}
        \item[(i)] $(r +1, c_1), (r+2, c_1), (r, c_3), (r+2,c_3) \in D$.
        \item[(ii)] $(r+2, c') \not\in D$ for all $c' > c_2$.
        \item[(iii)] $(r+1, c')\not\in D$ for all $c' > c_1$.
        \item[(iv)] $(r, c_1)\not\in D$.
    \end{itemize}
    \end{enumerate}
    Then $\mathcal{P}(D_0)$ is not ranked.
\end{lemma}

We can now prove our classification of northeast diagrams with ranked Kohnert posets. For the reader's convenience, we restate the theorem here. See Figure~\ref{fig:ranked_sub1} for a visualization of conditions (a) - (d) and Figure~\ref{fig:ranked_sub2} for an example of a northeast diagram whose Kohnert poset is not ranked. 

\begin{figure}[h]
\centering
\begin{subfigure}[t]{.55\textwidth}
  \centering
  \begin{tikzpicture}[scale=0.5]
    \celx(1,2)[x_1]
    \celx(1,4)[x_2]
    \celx(3,5)[x_3]
    \draw[pattern=dots] (1,0) rectangle (3,2);
    \draw[pattern=vertical lines] (3,0) rectangle (6,5);
    \draw (0,6)--(0,0)--(6,0)--(6,6);
    \node at (1.5,-0.5) {$c_1$};
    \node at (3.5,-0.5) {$c_3$};
    \node at (-0.5,2.5) {$r_1$};
    \node at (-0.5,4.5) {$r_2$};
    \node at (-0.5,5.5) {$r_3$};
    \end{tikzpicture}
  \caption{Figure~\ref*{fig:ranked_sub1}. The arrangement of cells $x_1$, $x_2,$ and $x_3$ satisfying (a) and (b). Condition (c) states that there is an empty position in each column of the dotted region. Condition (d) states that the number of cells in each column of the striped region is less than $r_1$. }
  \label{fig:ranked_sub1}
\end{subfigure}\qquad
\begin{subfigure}[t]{.4\textwidth}
  \centering
  \begin{tikzpicture}[scale = 0.5]
    \draw (1,4)--(1,0)--(3,0);
    \cel(1,1)
    \cel(1,2)
    \cel(2,2)
    \end{tikzpicture}
  \caption{Figure~\ref*{fig:ranked_sub2}. A minimal example of a northeast diagram whose Kohnert poset is not ranked.}
  \label{fig:ranked_sub2}
\end{subfigure}
\end{figure}

\begin{repeattheorem}[Theorem~\ref{thm:mainNEranked}]
     If $D_0$ is a northeast diagram, then $\mathcal{P}(D_0)$ is ranked if and only if $D_0$ does not contain $x_1 = (r_1,c_1), x_2 = (r_2,c_2)$, and $x_3 = (r_3,c_3)$ such that:
    \begin{itemize}
         \item[(a)] $r_1 < r_2 \leq r_3$,
         \item[(b)] $c_1 = c_2 < c_3$, 
         \item[(c)] for each  $c_1 \leq c <c_3$, there is at least one empty position $(r,c)$ where $r < r_1$, and
         \item[(d)] for each  $c \geq c_3$, the number of $r < r_3$ such that $(r, c) \in D_0$ is strictly less than $r_1$.
    \end{itemize}
 \end{repeattheorem}
 
\begin{proof}[Proof of Theorem~\ref{thm:mainNEranked}] (Proof of $\Rightarrow$) We show the contrapositive. Suppose $D_0$ contains cells $x_1$, $x_2$, and $x_3$ satisfying Properties (a) - (d). In the case that there are multiple such triples, order them lexicographically as tuples $(c_3, r_3, c_2, r_2, c_1, r_1)$ and choose the maximal tuple. This yields a choice of $x_1, x_2, x_3$ such that $x_3$ is right- and upper-most, then $x_2$ is the right- and upper-most choice of cell that works with $x_3$, then $x_1$ is the right- and upper-most choice of cell that works with $x_3$ and $x_2$.  
We now perform a sequence of Kohnert moves to obtain a diagram to which we may apply Lemma~\ref{lem:chmp}. This sequence is broken down into five subsequences, and we denote by $D_i$ the result after applying the $i$th subsequence of Kohnert moves.
\begin{enumerate}
    \item First we argue that there is at least one empty position $(r,c)$ with $r<r_1$ in every column $c$ with $c\geq c_3$. Supposing otherwise, by the northeast condition we would also have a cell at $(r_1,c)$, but this violates Property (d). Together with Property (c), this implies that we can perform Kohnert moves on $D_0$ to the cells strictly below and weakly to the right of $x_1$ to obtain a diagram $D_1$ that does not contain any cells in positions $(r_1 - 1, c)$ for $c \geq c_1$.
    
    \item By Property (d), we can perform Kohnert moves on $D_1$ to the cells strictly below and weakly to the right of $x_3$ to obtain a diagram $D_2$ that does not contain any cells in positions $(r, c)$ for $r_1 \leq r < r_3$ and $c \geq c_3$. Furthermore we observe that there are no cells to the right of $x_3$ in row $r_3$ of $D_2$. If such a cell $x_3'$ existed, the triple $(x_1, x_2, x_3')$ would certainly satisfy Properties (a), (b), and (d), and Property (d) applied to the triple $(x_1, x_2, x_3)$ would imply Property (c) for $(x_1, x_2, x_3')$. This would contradict the maximality of $(x_1, x_2, x_3)$. 
    
    \item By the maximality of our choices of $x_1,$ $x_2,$ and $x_3$, $D_2$ has at most one cell $(r, c)$ in each column $c_1 < c < c_3$ 
    with $r_1 \leq r \leq r_3$. By our previous set of moves, the only cell to the right of this rectangle is $x_3$. 
    Moreover, by our first set of moves, there are no cells in positions $(r_1-1,c)$ for any $c_1\leq c<c_3$. Then we can perform Kohnert moves on the cells in this rectangle of $D_2$ to obtain a diagram $D_3$ such that the only remaining cells $(r, c)$ with $r_1 \leq r < r_3$ and $c_1 \leq c$ are $x_1$ and $x_2$.
  \end{enumerate}
  
  We now have two cases, dependent on whether or not $(r_1-1,c_3)\in D_3$. We describe the rest of the subsequences and finish the proof separately for each case.
  
  \textbf{Case (a):}  $(r_1-1,c_3)\notin D_3$.
  \begin{enumerate}
   \item[(4a)] By the result of our second set of moves, we can perform $r_3 - r_1+1$ elementary moves on every cell $(r_3, c) \in D_3$ for $c_1< c \leq c_3$, 
   \emph{except for the leftmost such cell}, to which we apply $r_3-r_1$ moves. We refer to this leftmost cell, now in row $r_1$, as $x_3'=(r_1,c_3')$. If $(r_3, c) \in D_3$, then the maximality of our tuple implies that in the previous step there was no cell to move into position $(r_1 - 1, c)$. Thus, the moves described here are valid and yield a diagram $D_4$.

    \item[(5a)] As there are no cells $(r, c)$ in $D_4$
    with $r_1 < r \leq r_3$ and $c_1\leq c$ other than $x_2$, we can perform $r_2 - r_1 - 1$ elementary moves on $x_2$ to obtain a diagram $D_5$ such that the cell in column $c_1$ with label $x_2$ is in row $r_1 + 1$. We call this cell $x_2'$.
\end{enumerate}

We claim that $D_5$ satisfies the conditions (a)(i) through (a)(iv) of Lemma \ref{lem:chmp}(a), choosing cells $x_1, x_2'$, and $x_3'$ of $D_5$ (in that order) to be those described in condition (i) of the lemma. By our second set of moves, $D_5$ satisfies condition (ii) of the lemma. By our second, third, and final sets of moves, $D_5$ satisfies condition (iii) of the lemma. 
Finally, by our first set of moves and the hypothesis of this case, 
$D_5$ satisfies condition (iv) of the lemma. Then in this case, $\mathcal{P}(D_0)$ is not ranked.

\textbf{Case (b):} $(r_1-1,c_3)\in D_3$. 

If this cell has an empty space beneath it, we can again move all but the leftmost cell to the right of column $c_1$ in row $r_3$ down to row $r_1-1$, take this leftmost cell to the same $x_3'$ position, and apply Lemma~\ref{lem:chmp}(a) as above. So instead suppose there are no empty spaces beneath $(r_1-1,c_3).$

 \begin{enumerate}
   \item[(4b)] We perform $r_3 - r_1-1$ elementary moves on every cell $(r_3, c)$ for $c_1< c \leq c_3$, so that they all end up in row $r_1+1$. Call the resulting diagram $D_4$ and its cell $x_3'=(r_1+1,c_3)$.
   
    \item[(5b)] We perform $r_2 - r_1 - 1$ elementary moves on $x_2$ in $D_4$ to obtain a diagram $D_5$ such that the cell in column $c_1$ with label $r_2$ is in row $r_1 + 1$ as well. We call this cell $x_2'$.
\end{enumerate}
This diagram then satisfies the conditions (b)(i) through (b)(iv) in the statement of Lemma~\ref{lem:chmp}(b), choosing $x_1$, $x_2'$, $(r_1-1, c_3)$, and $x_3'$ of $D_5$ (in that order) to be the cells described in condition (b)(i). 
This completes the proof of the forward implication.

\hfill

(Proof of $\Leftarrow$) We show the contrapositive. Suppose $D_0$ is northeast and $\mathcal{P}(D_0)$ is not ranked. By Lemma \ref{lem:rankedCovering}, there are $D, D' \in \mathcal{P}(D_0)$ such that $D' \lessdot D$, but $D'$ is obtained from $D$ by a Kohnert move applied to a cell $x_2' = (r, c_1)$ such that $x_1' = (r-1, c_1) \in D$. If there were no cells of the form $x_3' = (r-1, c_3) \in D$ for some $c_3 > c_1$, then $D'$ could be obtained from $D$ by first applying a Kohnert move to $x_1'$ and then to $x_2'$. Thus, $D$ does contain such a cell $x_3'$.

Let $\mathcal{L}$ be the standard labeling of $D$ relative to $D_0$. For brevity, we say $L_1 = \mathcal{L}(x_1')$, $L_2 = \mathcal{L}(x_2')$, and $L_3 = \mathcal{L}(x_3')$. By strictness, $L_1 < L_2$. If $L_2 > L_3$, then by strictness and Property \ref{enum:cnortheast} of Definition \ref{def:NElabeling}, there is a cell $y = (r', c_3)$ in $D$ 
such that $r' > r - 1$ and $\mathcal{L}(y) = L_2$. By Property \ref{enum:dbelow} of Definition \ref{def:NElabeling}, $r' \leq r$. Thus, $r' = r$. In other words, $y$ is in the same row as $x_2'$ and to its right. This contradicts the fact that $x_2'$ is the rightmost cell in its row, which is necessary for it to be moved. Thus, $L_2 \leq L_3$.

By Theorem \ref{thm:labeling}, $D_0$ has cells $x_1 = (L_1, c_1)$, $x_2 = (L_2, c_1)$, and $x_3 = (L_3, c_3)$. From our construction of $c_1$ and $c_3$ and the derived conditions on $L_1, L_2$, and $L_3$, we have that $x_1$, $x_2$, and $x_3$ satisfy Properties (a) and (b) in the theorem statement. We now show that they satisfy Properties (c) and (d) and thus are the desired cells.

\begin{itemize}
    \item[(c)] A Kohnert move can be applied to $x_2' \in D$, so there is at least one empty space below $x_2' \in D$. Because $x_1'$ is one space below $x_2'$, that empty space is also below $x_1'$. By strictness and the column-equivalence of $\mathcal{L}$ to the super-standard labeling on $D_0$, $x_1 \in D_0$ and $x_1' \in D$ have the same number of cells below them in the same column. Then because $L_1 \geq r - 1$, $x_1$ must also have an empty space below it in the same column.

    Now fix a column $c$ with $c_1 < c < c_3$ and suppose that $(r', c) \in D_0$ for all $r'  < L_1$. 
    By the northeast property, $(L_1, c) \in D_0$ and $(L_2, c) \in D_0$. In particular, the cell $(L_1, c)$ has no empty spaces below it, so it can never move. Then $x_1$ cannot move either. In other words, $x_1 = x_1'$. By Property \ref{enum:dbelow} of Definition \ref{def:NElabeling}, the cell labeled $L_2$ in column $c$ of $D$ is weakly below $x_2'$, which is in the row directly above $x_1'$. But the only such available cell is to the right of $x_2'$, contradicting $x_2'$ being the rightmost cell in its row. Therefore, in each column $c$ with $c_1 < c < c_3$, there is at least one empty space in $D_0$ in a row below row $L_1$. 
    
    \item[(d)] For the sake of contradiction, suppose that this condition is not satisfied. Then there is a column $c \geq c_3$ 
    such that the total number of cells $(r' ,c)$ with $r' < L_3$ is at least $L_1$.
    By the northeast property, there is a cell $(L_3, c) \in D_0$. Then by strictness, the cell labeled $L_3$ in column $c$ of $D$ is above row $L_1$. Because $L_1 \geq r - 1$, the same cell is also strictly above row $r - 1$. But by Property \ref{enum:dbelow} of Definition \ref{def:NElabeling}, $x_3'$ would also need to be above row $r - 1$, contradicting the definition of $x_3'$ as $(r-1, c_3)$.
\end{itemize}
Then the initial diagram $D_0$ contains cells $x_1, x_2,$ and $x_3$ satisfying Properties (a) - (d) as described in the statement of the theorem.
\end{proof}

\section{Bounded posets}\label{sec:bounded}

\begin{figure}[h]
\centering
\begin{subfigure}[t]{.55\textwidth}
  \centering
  \begin{tikzpicture}[scale=0.5]
    \filldraw[color=blue!20] (1,0) rectangle (4,6);
    \celx(1,2)[x_1]
    \celx(4,3)[x_2]
    \celx(4,5)[x_3]
    \draw[pattern=dots] (1,0) rectangle (6,2);
    \draw[pattern=vertical lines] (1,3) rectangle (2,6);
    \draw (0,6)--(0,0)--(6,0)--(6,6);
    \node at (1.5,-0.5) {$c_1$};
    \node at (4.5,-0.5) {$c_2$};
    \node at (-0.5,2.5) {$r_1$};
    \node at (-0.5,3.5) {$r_2$};
    \node at (-0.5,5.5) {$r_3$};
    \end{tikzpicture}
  \caption{Figure~\ref*{fig:bounded_sub1}. The arrangement of cells $x_1$, $x_2,$ and $x_3$ satisfying (a) and (b). Condition (c) states that the number of cells in each column of the blue region is less than the number of cells in the column containing $x_2$ and $x_3$. Condition (d) states that there is an empty position in each column of the dotted region. Condition (e) states that there are no cells in the striped region.}
  \label{fig:bounded_sub1}
\end{subfigure}\qquad
\begin{subfigure}[t]{.4\textwidth}
  \centering
  \begin{tikzpicture}[scale = 0.5,baseline=-10pt]
    \draw (1,4)--(1,0)--(3,0);
    \cel(1,1)
    \cel(2,1)
    \cel(2,2)
    \end{tikzpicture}
  \caption{Figure~\ref*{fig:bounded_sub2}. A minimal example of a northeast diagram whose Kohnert poset is not bounded.}
  \label{fig:bounded_sub2}
\end{subfigure}
\end{figure}
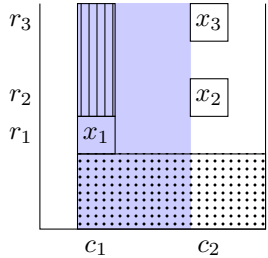

In this section, we prove Theorem~\ref{thm:mainNEbounded}, repeated here for the reader's convenience. We refer the reader to Figure \ref{fig:bounded_sub1} for a visualization of conditions (a) - (e) and to Figure~\ref{fig:bounded_sub2} for an example of a northeast diagram whose Kohnert poset is not bounded.

\begin{repeattheorem}[Theorem~\ref{thm:mainNEbounded}]
  If $D_0$ is a northeast diagram, then $\mathcal{P}(D_0)$ is bounded if and only if $D_0$ does not contain $x_1 = (r_1,c_1), x_2 = (r_2,c_2)$, and $x_3 = (r_3,c_3)$  such that:
    \begin{itemize}
        \item[(a)] $r_1 \leq r_2 < r_3$,
        \item[(b)] $c_1 < c_2 = c_3$,
        \item[(c)] for each $c_1 \leq c < c_2$, $\cwt(D_0)_c < \cwt(D_0)_{c_2}$,
        \item[(d)] for each $c\geq c_1$, there is at least one empty position $(r,c)$ where $r < r_1$, and 
        \item[(e)] for each $r_1 < r \leq r_3$, the cell $(r,c_1)$ is not in $D_0$.
    \end{itemize}
\end{repeattheorem}

By Proposition \ref{prop:NE elem chain no inversions}, every minimal element can be reached by a sequence of elementary moves with no move inversions, meaning that within the sequence, all elementary moves of $i$-cells occur before all elementary moves of $(i+1)$-cells for each $i \geq 1$. 
We claim that we can always construct a minimal element of $\mathcal{P}(D_0)$ for northeast diagram $D_0$ by applying the following ``greedy'' sequence of elementary moves with no move inversions to $D_0$. We assign the super-standard labeling to $D_0$. Then, for each $i$, starting with $i=2$, first move the rightmost $i$-cell downward as far as possible using only elementary moves. During this part of the process we call this $i$-cell the \textbf{moving cell} and say that this moving cell is taking its \textbf{turn}. Now, going leftward along the row, do the same for each subsequent $i$-cell. Then repeat the process for the $(i+1)$-cells. We refer to this greedy, inversion-free sequence of elementary moves as the \textbf{canonical procedure}.  In our next lemma, we prove our claim that this procedure produces a minimal element when $D_0$ is northeast. We also show that at no point during the canonical procedure is it ever possible to apply a jump Kohnert move to the cell that we are moving.

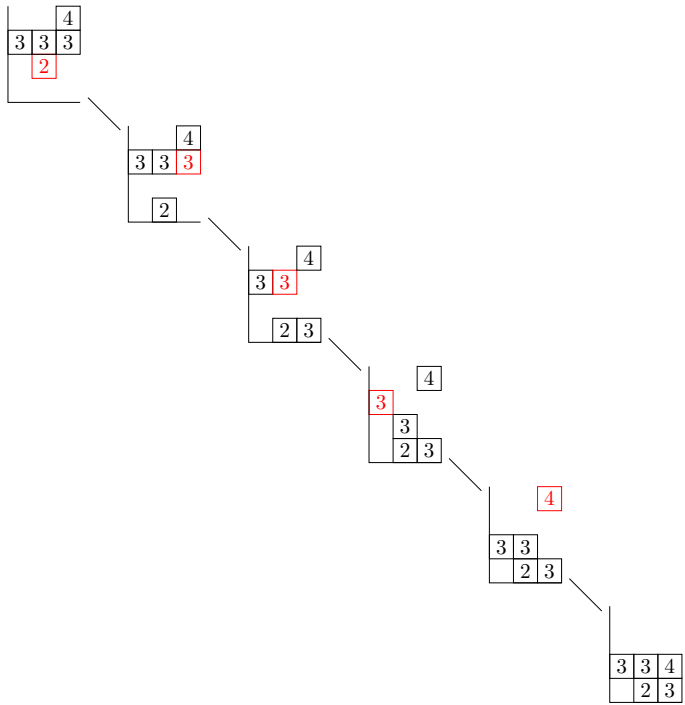
\begin{figure}[ht]
    \centering
    \scalebox{0.8}{\begin{tikzpicture}
 \node (1) at (-5,5) {\begin{tikzpicture}[scale=0.4]
    \rcelx(1,1)[2]
    \celx(0,2)[3]
    \celx(1,2)[3]
    \celx(2,2)[3]
    \celx(2,3)[4]
    \draw (0,4)--(0,0)--(3,0);
\end{tikzpicture}};
\node (2) at (-3,3) {\begin{tikzpicture}[scale=0.4]
    \celx(1,0)[2]
    \celx(0,2)[3]
    \celx(1,2)[3]
    \rcelx(2,2)[3]
    \celx(2,3)[4]
  \draw (0,4)--(0,0)--(3,0);
\end{tikzpicture}};
\node (3) at (-1,1) {\begin{tikzpicture}[scale=0.4]
    \celx(1,0)[2]
    \celx(0,2)[3]
    \rcelx(1,2)[3]
    \celx(2,0)[3]
    \celx(2,3)[4]
    \draw (0,4)--(0,0)--(3,0);
\end{tikzpicture}};
\node (4) at (1,-1) {\begin{tikzpicture}[scale=0.4]
    \celx(1,0)[2]
    \rcelx(0,2)[3]
    \celx(1,1)[3]
    \celx(2,0)[3]
    \celx(2,3)[4]
    \draw (0,4)--(0,0)--(3,0);
\end{tikzpicture}};
\node (5) at (3,-3) {\begin{tikzpicture}[scale=0.4]
    \celx(1,0)[2]
    \celx(0,1)[3]
    \celx(1,1)[3]
    \celx(2,0)[3]
    \rcelx(2,3)[4]
    \draw (0,4)--(0,0)--(3,0);
\end{tikzpicture}};
\node (6) at (5,-5) {\begin{tikzpicture}[scale=0.4]
    \celx(1,0)[2]
    \celx(0,1)[3]
    \celx(1,1)[3]
    \celx(2,0)[3]
    \celx(2,1)[4]
    \draw (0,4)--(0,0)--(3,0);
\end{tikzpicture}};
\draw (1)--(2)--(3)--(4)--(5)--(6);
\end{tikzpicture}}
    \caption{An example of constructing the canonical minimal element $D_{can}$ from an initial northeast diagram $D_0$ with the super-standard labeling. Note that the edges do not necessarily represent covering relations, for example in the last step.}\label{fig:canonical_diagram}
\end{figure}

For the remainder of this section, we let $D_{can}$ denote the element of $\mathcal{P}(D_0)$ that is obtained from the canonical procedure (see Figure \ref{fig:canonical_diagram} for an example). Throughout, $\mathcal{L}_{can}$ is the standard labeling of $D_{can}$ with respect to $D_0$. We refer to $D_{can}$ as the \textbf{canonical minimal element}, since it follows from Lemma \ref{lem:no jump and canonical minimal} that when $\mathcal{P}(D_0)$ is bounded, the diagram obtained in this way is its only minimal element. Note however that the minimality of $D_{can}$ is not immediate from the description of the canonical procedure and requires proof in Lemma \ref{lem:no jump and canonical minimal}. 

\begin{lemma}\label{lem:no jump and canonical minimal}
\!\!When applying the canonical procedure to a northeast diagram $D_0$, the following hold.
\begin{enumerate}
\item[(i)] At no point is it possible to apply a jump Kohnert move to the moving cell. In particular, once a cell ends its turn in the canonical procedure, it remains in that position in $D_{can}$.
\item[(ii)] The diagram $D_{can}$ resulting from the canonical procedure is a minimal element of $\mathcal{P}(D_0)$.
\end{enumerate}
\end{lemma} 

\begin{proof}
We begin with a proof of (i). Assume to the contrary that there is a moving cell to which a jump Kohnert move can be applied during the canonical procedure. Let $x$ be the first such cell. Let $(r,c)$ be the position of $x$ in some diagram, say $D'$, when its first jump Kohnert move is possible.
There must be an empty position $(r',c)$ where $r' < r-1$, and there must be cells $x_1, \dots, x_n$ that occupy all positions $(r'+1, c), \dots, (r-1, c)$, respectively (with $n=r-r'-1$). Further, $x$ must be the rightmost cell in row $r$ of $D'$.

For each $x_i$ with $1 \leq i \leq n$, consider the diagram $D_i$ obtained immediately after $x_i$ has ended its turn as the moving cell. By the minimality of $x$, a jump Kohnert move could not have been applied to move $x_i$ into the position $(r'+i,c)$. Thus, it must have been the case that $x_i$ had a cell to its right in the same row in $D_i$. Let $y_i$ be the leftmost such cell in $D_i$. Let $\mathcal{L}_i$ be the standard labeling of $D_i$ with respect to $D_0$. Then Lemma~\ref{lemma:above} gives $\mathcal{L}_i(y_i) \geq \mathcal{L}_i(x_i)$, while by the definition of the canonical procedure, $\mathcal{L}_i(y_i) \leq \mathcal{L}_i(x_i)$ as we have not yet moved any cells above row $\mathcal{L}_i(x_i)$. Thus, $\mathcal{L}_i(y_i) = \mathcal{L}_i(x_i)$ for each $i$.
 
In particular, $\mathcal{L}_n(y_n) = \mathcal{L}_n(x_n)< \mathcal{L}_{n}(x)$.
Thus, by Property (3) of Definition \ref{def:NElabeling}, there must be a cell $y$ in the same column as $y_{n}$ in $D'$ with $\mathcal{L}_{n}(y) = \mathcal{L}_{n}(x)$. By Property (4) of Definition \ref{def:NElabeling},
$y$ must be located in row $r$ to the right of $x$ in $D'$. This contradicts our assumption that $x$ can move and completes our proof of (i).

Now we show that $D_{can}$ is minimal. We know that $D_{can}$ is a minimal element of $\mathcal{P}(D_0)$ if and only if no cells in $D_{can}$ are movable. By part (i), during the canonical procedure, each moving cell ends its turn immovable.
This implies that each moving cell ends its turn with either a cell to the right of it in its row or with no empty spaces beneath it in its column.

Suppose there exists a cell $x \in D_{can}$ that is movable. Then there is an empty space somewhere below $x$. By (i), this empty space must have existed when $x$ ended its turn as the moving cell, yielding a diagram $D'$. Thus, there must have been some cell $y$ in the same row to the right of $x$ in $D'$. 
But then by (i), the cell $y$ must be in the same position in $D_{can}$ as it is in $D'$. This contradicts that $x$ is movable and completes our proof of (ii).
\end{proof}

Since $\mathcal{P}(D_0)$ is not bounded if and only if it has two or more distinct minimal elements, to prove Theorem \ref{thm:mainNEbounded}, it is sufficient to prove the following proposition.

\begin{proposition}\label{cor:mainNEboundedredux}
Let $D_0$ be a northeast diagram. Then $D_0$ contains cells $x_1, x_2, x_3$ satisfying conditions (a) - (e) of Theorem~\ref{thm:mainNEbounded} if and only if $\mathcal{P}(D_0)$ contains a minimal element $\tilde{D}$ distinct from the canonical minimal element $D_{can}$.  
\end{proposition}

The proof of Proposition \ref{cor:mainNEboundedredux} requires a few additional lemmas. Our first is a general result about minimal diagrams.

\begin{lemma}\label{lem:bdd claim 2}
Let $\tilde{D}$ be a minimal element of $\mathcal{P}(D_0)$. Suppose there exists a column $c$ that contains an empty position below a cell. 
Let $r_{\max}$ be the maximum row index such that $(r_{\max},c) \in \tilde{D}$. Then there exists a column to the right of column $c$ in $\tilde{D}$ in which all rows $1$ through $r_{\max}$ are occupied by cells.
\end{lemma} 

\begin{proof}
Since $\tilde{D}$ is a minimal element, there must be at least one column $c' > c$ such that $(r_{\max},c') \in \tilde{D}$. Let $c_{\max}$ equal the maximum such column. As $(r_{\max}, c_{\max})$ is the rightmost cell in its row, the minimality of $\tilde{D}$ implies that in column $c_{\max}$, all rows $1$ through $r_{\max}-1$ are occupied by cells. This completes our proof, since $(r_{\max}, c_{\max}) \in \tilde{D}$. 
\end{proof}

The next lemma specializes Lemma~\ref{lem:bdd claim 2} to the canonical minimal element to deduce further information about $D_0$. We introduce the notation $\cwt_{\leq r}(D)_c$ (respectively, $\cwt_{< r}(D)_c$) to denote the number of cells in column $c$ of diagram $D$ that lie weakly (respectively, strictly) below row $r$. Likewise, let $\cwt_{\geq r}(D)_c$ (respectively, $\cwt_{> r}(D)_c$) denote the number of cells in column $c$ of diagram $D$ that lie weakly (respectively, strictly) above row $r$.

\begin{lemma} \label{lem:space-below-cell-iff-heavy-column-at-right} Let $x=(r,c)$ be a cell in the northeast diagram $D_0$ that reaches position $x'=(r',c)$ at the end of the canonical procedure in the canonical minimal element $D_{can}$. If $x'$ has an empty position below it in column $c$ of $D_{can}$, then:
\begin{enumerate}
\item[(i)] there is a column ${c'}> c$ that satisfies $\cwt_{\leq r}(D_0)_{{c'}}>\cwt_{\leq r}(D_0)_{c}$ and
\item[(ii)] if $\hat{c}$ is the leftmost such column, then the cell $\hat{x} = (r', \hat{c})$ exists in $D_{can}$ and $\mathcal{L}_{can} (\hat{x})=r$. 
\end{enumerate}
\end{lemma}

\begin{proof} 
First we show (i). Let $D'$ be the diagram reached after $x$ has ended its turn as the moving cell in position $(r',c)$. Since $x'$ has an empty position below it in column $c$ of $D_{can}$, it must have an empty position below it in column $c$ of $D'$ by part (i) of Lemma \ref{lem:no jump and canonical minimal}. This, combined with the fact that $x'$ is not movable by part (i) of Lemma \ref{lem:no jump and canonical minimal} implies that $x'$ must have a cell in the same row to the right of it in $D'$. Let $y'$ be the rightmost such cell and let $\mathcal{L}'$ be the standard labeling on $D'$ with respect to $D_0$. Then Lemma~\ref{lemma:above} gives $\mathcal{L}'(y')\geq \mathcal{L}'(x')=r$, while by the definition of the canonical procedure, $\mathcal{L}'(y')\leq r$ as we have not yet moved any cells above row $r$. Hence, $y'$ is an $r$-cell.

Suppose $y'$ is in column ${c}'$. There must not be any empty spaces below $y'$ since $y'$ is not movable and it is the rightmost cell in its row in $D'$. Thus, $D'$ has cells labeled by exactly $r'$ values less than or equal to $r$ in column ${c}'$ while $D'$ has strictly fewer than $r'$ cells labeled less than or equal to $r$ in column $c$. We conclude that $\cwt_{\leq r}(D_0)_{{c}'}>\cwt_{\leq r}(D_0)_{c}$, which proves (i).

Next, we show (ii). Let $\hat{c}$ be the leftmost column satisfying $\cwt_{\leq r}(D_0)_{\hat{c}} > \cwt_{\leq r}(D_0)_{c}$ such that $\hat{c}>c$. From (i), we know that such a column exists, and is either $c'$ or to the left of $c'$, so $c<\hat{c}\leq c'$. Because of the column-weight condition, there is at least one cell in column $\hat{c}$ in a row strictly lower than $r$ in $D_0$. So, since $D_0$ is northeast, $(r,\hat{c}) \in D_0$, or equivalently, there is an $r$-cell in column $\hat{c}$. 
Then, in $D_{can}$, the equality $\mathcal{L}_{can}(x')=r=\mathcal{L}_{can}(y')$ forces the $r$-cell in column $\hat{c}$ to be at $(r',\hat{c})$ by Property (4) of Definition \ref{def:NElabeling}.  Thus, $\hat{x}=(r',\hat{c})\in D_{can}$ and $\mathcal{L}_{can} (\hat{x})=r$, completing the proof of (ii). 
\end{proof}

Our next lemma consists of the most technical parts of the proof for the $(\Rightarrow)$ direction of Proposition~\ref{cor:mainNEboundedredux}. 

\begin{lemma}\label{lem:bdd claim 1}
Let $x_1$, $x_2$, and $x_3$ satisfy conditions (a) - (e) of Theorem~\ref{thm:mainNEbounded} and be chosen so that the tuple $(c_1, r_1, -c_2, r_2, -c_3, r_3)$ is maximal in the lexicographic order on $\Z^6$. Indicate the positions of the $x_i$ in $D_{can}$ at the end of the canonical procedure by $x_i'=(r_i',c_i)$. Then:
\begin{enumerate}[label=(\roman*)]
\item for $r>r_3$, $(r,c_1) \in D_0$ if and only if $(r,c_2) \in D_0$, and 
\item $r_1'\leq r_2'$.
\end{enumerate}
\end{lemma}

\begin{proof}
Notice that the choice of $x_1$, $x_2$, and $x_3$ forces $x_1$ to have the maximum possible $c_1$, then maximum possible $r_1$. With respect to this $x_1$, this then forces the minimum possible $c_2$, and then the maximum possible $r_2$. Finally, with respect to the chosen $x_1$ and $x_2$, this forces the the maximum possible $r_3$.

We start with part (i). If $r>r_3$ and $(r,c_1) \in D_0$, then since $D_0$ is northeast, $x_2$ and $x_3$ force $(r,c_2) \in D_0$. For the converse, let $y=(r_y,c_2)$ be a cell with $r_y>r_3$. If $(r_y, c_1)\not\in D_0$, then there are two cases. If $(r,c_1) \not\in D_0$ for all $r_3<r<r_y$, then Properties (a) - (e) hold for $x_1,x_2,y$. As $y = (r_y,c_2)$ has $r_y>r_3$, this contradicts the maximality of the chosen $x_1$, $x_2$, and $x_3$. If there does exist $(r,c_1)\in D_0$ with $r_3<r<r_y$, then let $z=(r_z,c_1)$ be the upper-most such cell. Since $D_0$ is northeast, $x_2$ and $x_3$ force $(r_z,c_2)\in D_0$. Then $z$, $(r_z,c_2)$, and $y$ form a triple that satisfies (a) - (e), again contradicting the maximality of the choice $x_1$, $x_2$, $x_3$. So we must have $(r_y,c_1) \in D_0$, and thus part (i) follows.

We now prove part (ii), continuing with $x_1$, $x_2$, and $x_3$ chosen so that $(c_1,r_1,-c_2,r_2,-c_3,r_3)$ is maximal in the lexicographic ordering. Part (i) and condition (e) imply that all of the cells in column $c_2$ above row $r_3$ have a corresponding cell in column $c_1$ above row $r_3$. This means that $x_3$ is the only cell in column $c_2$ above row $r_2$ that does not have a cell in the same row in column $c_1$.

In other words, $\cwt(D_0)_{c_1}- \cwt_{\leq r_1}(D_0)_{c_1}+1=\cwt(D_0)_{c_2}-\cwt_{\leq r_2}(D_0)_{c_2}$.
Rewriting this yields \[
(\cwt_{\leq r_2}(D_0)_{c_2}+1) - \cwt_{\leq r_1}(D_0)_{c_1} = \cwt(D_0)_{c_2} - \cwt(D_0)_{c_1} > 0,
\] where the fact that this quantity is greater than zero is due to condition (c). Thus, we have that \begin{equation} \label{eq:bwt}\cwt_{\leq r_1}(D_0)_{c_1}\leq \cwt_{\leq r_2}(D_0)_{c_2}. \end{equation}
If $x_1'$ has no empty positions beneath it, then $r_1'\leq r_2'$ follows directly from  \eqref{eq:bwt}. Hence we may now assume $x_1'$ has a space beneath it in its column. By Lemma~\ref{lem:space-below-cell-iff-heavy-column-at-right}, this means that there is some column $\hat{c}>c_1$ with 
\begin{equation}
\label{eq:hatcdef}
\cwt_{\leq r_1}(D_0)_{\hat{c}}> \cwt_{\leq r_1}(D_0)_{c_1}. 
\end{equation}

We choose the leftmost such column, and
we claim that such a $\hat{c}$ must satisfy $\hat{c}\geq c_2$. Suppose, for contradiction, that $\hat{c}$ satisfies $c_1<\hat{c}<c_2$. 

\textbf{Case 1:} $(r_3, \hat{c})\notin D_0$. Let $z=(r_z,\hat{c})$ be the uppermost cell in column $\hat{c}$ that is below $r_3$; such a cell must exist by  \eqref{eq:hatcdef}. Then either the triple $z$, $x_2$, $x_3$ satisfies (a) - (e) if $r_z\leq r_2$, or by the northeast property of $D_0$ there is a triple $z$, $(r_z,c_2)$, $x_3$ that satisfies (a) - (e). In either situation, lexicographic maximality of our choice of triple is violated since $\hat{c}>c_1$, which produces a contradiction.

\textbf{Case 2:} $(r_3, \hat{c})\in D_0$. 
Since $\cwt_{\leq r_1}(D_0)_{\hat{c}}> \cwt_{\leq r_1}(D_0)_{c_1}>0$, by the northeast property of $D_0$ we have that $(r_1,\hat{c})\in D_0$ so there is at least one cell of the form $(s, \hat{c})$ with $r_1\leq s <r_3$. Let $y=(r_y, \hat{c})$ be the cell of this form with maximal row index $r_y$. Then the triple $x_1$, $y$, $(r_3,\hat{c})$ straightforwardly satisfies conditions (a), (b), (d), and (e). We now show that condition (c) is satisfied as well. 

The combination of $(r_1,\hat{c})\in D_0$ and the northeast property of $D_0$ imply that if $(s, c_1) \in D_0$, then $(s, \hat{c}) \in D_0$ for all $s > r_1$. Thus, $\cwt_{> r_1}(D_0)_{\hat{c}}\geq \cwt_{> r_1}(D_0)_{c_1}$. Adding this inequality to \eqref{eq:hatcdef} gives
\begin{equation}\label{eq:condition-c-part1}
\cwt(D_0)_{c_1}<\cwt(D_0)_{\hat{c}}.
\end{equation}
To handle the other columns, suppose for sake of contradiction that the inequality 
\begin{equation}\label{eq:part-c-for-contradiction}
\cwt(D_0)_{c'}=\cwt_{\leq r_1}(D_0)_{c'}+\cwt_{> r_1}(D_0)_{c'}\geq \cwt_{\leq r_1}(D_0)_{\hat{c}}+\cwt_{> r_1}(D_0)_{\hat{c}}=\cwt(D_0)_{\hat{c}},
\end{equation}
holds for some $c_1<c'<\hat{c}$. 
Once again, the combination of $(r_1,\hat{c})\in D_0$ and the northeast property of $D_0$ imply
\[\cwt_{> r_1}(D_0)_{c'}\leq\cwt_{> r_1}(D_0)_{\hat{c}} \]
Combined with \eqref{eq:part-c-for-contradiction}, this implies that 
\[\cwt_{\leq r_1}(D_0)_{c'}\geq\cwt_{\leq r_1}(D_0)_{\hat{c}}.
\]
But since $c'<\hat{c}$, this contradicts our choice of $\hat{c}$ as the least column index $c>c_1$ satisfying $\cwt_{\leq r_1}(D_0)_{c}\geq \cwt_{\leq r_1}(D_0)_{c_1}$. We must then conclude that $\cwt(D_0)_{c'}<\cwt(D_0)_{\hat{c}}$ for each $c'$ such that $c_1<c'<\hat{c}$. Together with \eqref{eq:condition-c-part1}, this shows that condition (c) applies to the triple $x_1$, $y$, $(r_3, \hat{c})$ as well, contradicting lexicographic maximality of our original triple and completing this case.

Having dealt with both cases, we have now shown that $\hat{c}\geq c_2$. Applying Lemma~\ref{lem:space-below-cell-iff-heavy-column-at-right}(ii), we know that there is an $r_1$-cell $(r_1',\hat{c})\in D_{can}$. If $\hat{c}=c_2$, then the fact that $r_2'\geq r_1'$ follows since $x_2$ starts above this $r_1$-cell in $D_0$. Otherwise if $\hat{c}>c_2$, the fact $r_2'\geq r_1'$ follows from Lemma~\ref{lemma:above}, comparing the labels on $x_2'$ and $(r_1', \hat{c})$. This completes the proof of (ii). 
\end{proof}

We are now ready to prove Theorem \ref{thm:mainNEbounded} by proving Proposition \ref{cor:mainNEboundedredux}.

\begin{proof}[{Proof of Theorem~\ref{thm:mainNEbounded}}]
We show both $\Rightarrow$ and $\Leftarrow$ directions of Proposition~\ref{cor:mainNEboundedredux}.

(Proof of $\Rightarrow$) Assume $D_0$ contains cells $x_1, x_2, x_3$ satisfying conditions (a) - (e). As in Lemma~\ref{lem:bdd claim 1}, we choose the triple such that $(c_1, r_1, -c_2, r_2, -c_3, r_3)$ is maximal in the lexicographic order. Assign to $D_0$ the super-standard labeling, and to each diagram in $\mathcal{P}(D_0)$ the standard labeling with respect to $D_0$.  We recall the notation $x_i'=(r_i',c_i)$ for the cell in column $c_i$ of $D_{can}$ with label $r_i$. We will construct an element of $\mathcal{P}(D_0)$ called $\tilde{D}$, and show that $\tilde{D}$ is a minimal element of $\mathcal{P}(D_0)$ distinct from the canonical minimal element $D_{can}$. We let $\tilde{x}_i = (\tilde{r}_i,c_i)$ denote the cell in $\tilde{D}$ that originated as $x_i$ in $D_0$. 

We first show that during the canonical procedure, $x_1$ moves down at least one position. Assume to the contrary that it does not move,
and recall that by condition (d), each column $c \geq c_1$ has at least one empty position $(r,c)$ for $r<r_1$ in $D_0$.
Thus, there is an empty position somewhere below $x_1$ in $D_0$ and also in the diagram $D'$ reached at the start of $x_1$'s turn as the moving cell. By Lemma~\ref{lem:no jump and canonical minimal}, if a cell does not move in the canonical procedure, it is not because there is only a jump move available. So, in $D'$, the cell $x_1$ is not rightmost in its row. Let $y_1$ be the rightmost cell in row $r_1$ of $D'$. Let $\mathcal{L}'$ be the standard labeling of $D'$. Then Lemma~\ref{lemma:above} gives $\mathcal{L}'(y_1)\geq \mathcal{L}'(x_1)=r_1$, while by the definition of the canonical procedure, $\mathcal{L}'(y_1)\leq r_1$ as we have not yet moved any cells above row $r_1$. Hence, $y_1$ is an $r_1$-cell. However, since $y_1$ is in row $r_1$, this means that $y_1$ did not move in the canonical procedure, and yet has an empty position below it. Again by Lemma~\ref{lem:no jump and canonical minimal}, at no point is it possible to apply a jump Kohnert move to a cell, so there is an empty position immediately below $y_1$. This contradicts that $y_1$ is the rightmost cell in its row, and therefore $x_1$ must have moved at least once during the canonical procedure.

To construct $\tilde{D}$, we begin applying elementary moves to $D_0$ according to the canonical procedure, but refrain from making the final possible elementary move on the cell in column $c_1$ with label $r_1$, so that it instead rests in position $(r_1'+1,c_1)$, leaving an empty space in position $(r_1',c_1)$. Denote the diagram obtained immediately after this cell with label $r_1$ stops moving at row $r_1'+1$ as $D''$. From $D''$, we continue to move all subsequent cells greedily down as far as possible using only elementary moves, continuing from the cell immediately left of $x_1$ as in the canonical procedure. At the end of this procedure, if necessary, in column $c_1$ apply Kohnert moves to the cells strictly above row $r_1'+1$ to move these cells down as far as possible until they are no longer movable. We will see that none of these cells can jump below $(r_1'+1,c_1)$, and thus  position $(r_1',c_1)$ remains empty. Then, if any cells in columns strictly to the left of column $c_1$ are movable, apply Kohnert moves in these columns until this is no longer the case. After all such moves are made, we call the resulting diagram $\tilde{D}$ and denote the position of the $r_i$-cell in column $c_i$ (the cell that started as $x_i$) by $\tilde{x}_i = (\tilde{r}_i,c_i) \in \tilde{D}$ for $i\in\{1,2,3\}$. 
It remains to verify that $D_{can}$ and $\tilde{D}$ differ in column $c_1$, and that $\tilde{D}$ is minimal.

First we show that $D_{can}$ and $\tilde{D}$ differ in column $c_1$ by showing that no cell above $(r_1'+1,c_1)$ in $D''$ can jump over $(r_1'+1,c_1)$ into the empty position $(r_1',c_1)$.
Each column strictly to the right of column $c_1$ is identical in $D_{can}$ and $\tilde{D}$, as these columns are unaffected by changes in column $c_1$. We use this fact repeatedly in what follows. By Lemma~\ref{lem:bdd claim 1}, $r_1' \leq r_2'$, so 
\begin{equation} 
\label{eq:key ineq bdd}
r'_1 +1\leq r'_2+1 = \tilde{r}_2 + 1 \leq  \tilde{r}_3 = r'_3.
\end{equation}
Suppose there is another cell above $x_1$ in $D_0$, and let $y_1=(r_y,c_1)$ be such a cell with minimal $r_y$, noting that $r_y>r_3$ by condition (e). By Lemma~\ref{lem:bdd claim 1}(i), there is a corresponding $y_2=(r_y,c_2)\in D_0$. In $\tilde{D}$ (and $D_{can}$), $y_2$ must land at some row ${r}'_{y}> r_3'$. But then by Lemma~\ref{lemma:above}, $y_1$ must land in a row with index at least ${r}'_y$, which is then strictly greater than $r_1'+1$ by \eqref{eq:key ineq bdd}.

Thus, in applying additional Kohnert moves to the cells in $D''$ that lie above $(r_1'+1,c_1)$, none of these cells land in position $(r_1',c_1)$. Therefore, position $(r_1',c_1)$ is occupied in $D_{can}$ but empty in $\tilde{D}$. This proves that $D_{can}$ and $\tilde{D}$ differ in column $c_1$. 

Next we show that $\tilde{D}$ is minimal. Our construction of $\tilde{D}$ ensures that no cells in column $c_1$ strictly above $\tilde{x}_1$ or in any column strictly to the left of column $c_1$ are movable. Furthermore, the fact that $\tilde{D}$ and $D_{can}$ are identical in columns to the right of column $c_1$ as well as below row $r_1'$ of column $c_1$ means that none of those cells are movable either. Thus, to show that $\tilde{D}$ is minimal, we only need to show that $\tilde{x}_1$ is not the rightmost cell in its row.

If all rows $1 \leq r \leq \tilde{r}_3$ of column $c_2$ are occupied in $\tilde{D}$, then $\tilde{x}_1$ is not rightmost in its row and is immovable. Otherwise, the fact that $\tilde{x}_3$ has a space somewhere below it in $\tilde{D}$ implies that $x'_3$ has a space somewhere below it in $D_{can}$. Lemma \ref{lem:bdd claim 2} then implies that there is some column $c'>c_2$ in $D_{can}$ that has rows $1$ through $r'_3$ occupied by cells. Thus, $\tilde{D}$ also has rows $1$ through $\tilde{r}_3 = r'_3$ occupied by cells in column $c'$. Hence, in this case it is also true that $\tilde{x}_1$ is not rightmost in its row.  This completes our proof of ($\Rightarrow$).

\hfill

(Proof of $\Leftarrow$) Assume that $\mathcal{P}(D_0)$ has a second minimal element, $\tilde{D}$, distinct from $D_{can}$. This $\tilde{D}$ can be produced by some sequence of elementary Kohnert moves with no move inversions by Proposition \ref{prop:NE elem chain no inversions}. We show that there exist cells $x_1, x_2, x_3$ in $D_0$ satisfying conditions (a) - (e). 

Let $\mathcal{L}_{can}$ and $\tilde{\mathcal{L}}$ denote the standard labelings on $D_{can}$ and $\tilde{D}$, respectively. Let column $c_1$ be the rightmost column where $D_{can}$ and $\tilde{D}$ differ. Let $s_1'$ be the smallest row index where $D_{can}$ and $\tilde{D}$ differ in column $c_1$.

We claim that the position $(s_1', c_1)$ is empty in $\tilde{D}$ and occupied in $D_{can}$. For sake of contradiction, suppose otherwise that $x=(s_1',c_1)\in\tilde{D}$ and is empty in $D_{can}$. Then, since both diagrams have the same column weights, there is a cell $x'=(r',c_1)$ in $D_{can}$ that is the first occupied space above row $s_1'$. Since $x'$ has an empty position below it, by Lemma~\ref{lem:space-below-cell-iff-heavy-column-at-right}(ii), there is a column $\hat{c}$ to the right of column $c_1$ and a cell $\hat{x}=(r',\hat{c})\in D_{can}$ with $\mathcal{L}_{can}(\hat{x})=\mathcal{L}_{can}(x')$. We also have $\tilde{\mathcal{L}}(x)=\mathcal{L}_{can}(x')$, and since $D_{can}$ and $\tilde{D}$ are identical to the right of column $c_1$, $\tilde{\mathcal{L}}(\hat{x})=\mathcal{L}_{can}(\hat{x})$ as well. But this implies that $x$ and $\hat{x}$ are cells of $\tilde{D}$ with the same label and for which $\hat{x}$ is above and to the right of $x$, contradicting Lemma~\ref{lemma:above}. Thus, $(s_1',c_1)$ must be occupied in $D_{can}$ and empty in $\tilde{D}$.

Let $s_1 = \mathcal{L}_{can}((s_1',c_1))$. Then the first cell in column $c_1$ of $\tilde{D}$ above row $s_1'$ is precisely the cell with label $s_1$. Call this cell $\tilde{y}_1=(\tilde{s}_1,c_1)$, let $y_1 = (s_1,c_1)\in D_0$, and denote by $y_1'=(s_1',c_1)$ the $s_1$-cell in $D_{can}$.

For each nonnegative integer $j$, let $C_j:=\{c_1,c_{1}+1,\dots,c_1+j\}$ be a set of consecutive column indices. Let $(r_{C_j}, c_{C_j})$ be the lexicographically maximum cell in $\tilde{D}$ such that $c_{C_j} \in C_j$ and the cell has an empty position somewhere below it; this maximum exists since $\tilde{y}_1$ is in column $c_1$ and has an empty position below it. Notice that by this lexicographic maximality, as $j$ increases, the position of the cell $(r_{C_j},c_{C_j})$ moves weakly rightward and upward. For each $j$, by applying Lemma~\ref{lem:bdd claim 2}, there is a leftmost column $m_j > c_{C_j}$ in which all rows 1 through $r_{C_j}$ contain cells. Then as $j$ increases, $m_j$ weakly increases. Let $j'$ be the minimum index such that $m_{j'} = c_1+{j'}$. Such an index always exists since the sequence $(m_0, m_1, m_2,\dots)$ is weakly increasing and bounded above, with $c_1<m_j\leq m$ for all $j$, where $m$ is the last column. Further let $c_2 = m_{j'}$ and $r_{\max} = r_{C_{j'}}$.

We point out a useful property of this choice of $c_2$. Notice that column $c_2$ of $\tilde{D}$ contains cells in rows 1 through $r_{\max}$. Furthermore, for each column $c$ such that $c_1 \leq c < c_2$, we claim that the maximum row index of any cell in column $c$ of $\tilde{D}$ is also $r_{\max}$. To show this, suppose for contradiction that there is a cell $(r,c)\in \tilde{D}$ with $r>r_{\max}$ and $c_1\leq c <c_2$. If there are no empty positions below $(r,c)$, then consider columns $c_1$ through $c=c_1+d$. We have that $m_0,\ldots, m_d$ are each at most $c$, and since this sequence is weakly increasing, there exists some $m_{j''} = c_1+j''$ further left than $j'$, which contradicts the definition of $j'$. If there is an empty position below $(r,c)$, then $(r,c)$ would be a lexicographically larger cell than $(r_{C_{j'}},c_{C_{j'}})$, contradicting our choice when identifying this cell. 

Therefore, all cells in each column $c_1 \leq c < c_2$ are contained in rows 1 through $r_{\max}$. By a similar argument to the previous paragraph, if in one of these columns, every such row contains a cell, we would contradict the choice of $m_{j'}$. Thus,\begin{equation} \label{eq:rmax column weights}
\cwt_{\leq r_{\max}}(\tilde{D})_c < r_{\max} = \cwt_{\leq r_{\max}}(\tilde{D}_{c_2})\,\,\mathrm{for\,\, all}\,\,c_1\leq c < c_2.
\end{equation}

Since $(r_{\max},c_{C_{j'}})$ is a cell with a space beneath it, $r_{\max} \geq 2$. Hence, column $c_2$ contains at least two cells.
We claim that $D_0$ contains two cells $y_2=(s_2,c_2)$ and $y_3=(s_3,c_3)$ where $s_1 \leq s_2 <s_3$ and $c_1 < c_2 = c_3$. As the canonical procedure that yields $D_{can}$ moves $y_1$ down at least one row below the position of $\tilde{y}_1$ in $\tilde{D}$, the cell $(\tilde{s}_1-1,c_2)$ in  $\tilde{D}$ has label greater than or equal to $\tilde{\mathcal{L}}(\tilde{y}_1) = s_1$ by Lemma~\ref{lemma:above}. Call this cell $\tilde{y}_2$ and let $s_2 = \tilde{\mathcal{L}}(\tilde{y}_2)$ Then by the strictness on the labels in each column, the cell directly above $\tilde{y}_2$, which we call $\tilde{y}_3$, has $s_3 = \tilde{\mathcal{L}}(\tilde{y}_3)$ where $s_3>s_2 \geq s_1$.  As these labels are equal to the initial row indices of the corresponding cells in $D_0$ (which we refer to as $y_2$ and $y_3$, respectively), $D_0$ contains cells $y_2=(s_2,c_2)$, $y_3=(s_3,c_2)$ with $s_1 \leq s_2 <s_3$ and $c_1 < c_2 = c_3$, proving our claim. 

Having shown that there exists in $D_0$ a triple of cells $y_1$, $y_2$, $y_3$ that satisfies (a), (b), and the inequality $\tilde{s}_1>s_1'$, we now choose a particular such triple in the same columns.  Assign $x_1= (r_1, c_1)$, $x_2=(r_2,c_2)$, and $x_3=(r_3,c_3)$ to be the cells in columns $c_1$ and $c_2=c_3$ of $D_0$ such that the tuple $(r_1,-r_2,-r_3)$ is lexicographically maximal, conditions (a) and (b) are satisfied, and $\tilde{r}_1>r_1'$, where $\tilde{r}_1$ and $r_1'$ are the row indices of the $r_1$-cells in column $c_1$ of $\tilde{D}$ and $D_{can}$, respectively. 
We now show that this choice of $x_1,$ $x_2,$ and $x_3$ either satisfies conditions (a) - (e) or can be used to produce a related triple that does, particularly in the event that (e) is not satisfied.

\begin{itemize}
\item[(a),(b)] These inequalities are satisfied by $x_1$, $x_2,$ and $x_3$ by construction.
\item[(c)] 
We saw in (\ref{eq:rmax column weights}) that $\cwt_{\leq {r_{\max}}}(\tilde{D})_c < r_{\max} = \cwt_{\leq {r_{\max}}}(\tilde{D})_{c_2}$ for all $c_1 \leq c < c_2$. 
Further, we have seen that all cells in $\tilde{D}$ in columns $c$ with $c_1 \leq c < c_2$ are located in rows $1$ through $r_{\max}$. 
So $\cwt_{\leq {r_{\max}}}(\tilde{D})_c = \cwt(D_0)_c$ for each $c_1 \leq c < c_2$. 
Thus, $$\cwt(D_0)_c < \cwt_{\leq {r_{\max}}}(\tilde{D})_{c_2} \leq \cwt(D_0)_{c_2}$$ for all $c_1 \leq c < c_2$, so condition (c) is satisfied. 

\item [(d)] For both $D_{can}$ and $\tilde{D}$ to be reachable from $D_0$, by Lemma~\ref{lemma:above}, it must be the case that all $r_1$-cells in columns $c \geq c_1$ in $D_0$ move down during the canonical procedure. 
This means that each column $c \geq c_1$ that contains an $r_1$-cell must have an empty position in some row $r < r_1$ in $D_0$. Otherwise, if a column $c \geq c_1$ does not contain an $r_1$-cell, then assume for the sake of contradiction that there is no empty position $(r,c)$ in $D_0$ for some $r < r_1$. Then as every $(r,c)$ with $r<r_1$ is a cell in $D_0$, the presence of $x_1$ implies that $(r_1,c) \in D_0$ by the northeast property of $D_0$. This contradicts our assumption that the column had no $r_1$-cell and hence this case can not occur. Thus, condition (d) is satisfied. 

\item[(e)] Finally, we consider cells in positions $(r, c_1)$ in $D_0$ such that $r > r_1$. We must show that $D_0$ does not contain any cells $(r,c_1)$ with $r_1 < r \leq r_3$. If there was a cell $(r,c_1) \in D_0$ with $r_2 < r < r_3$, then the northeast property of $D_0$ would imply the presence of a cell in column $c_2$ that would contradict our choice of $x_3$. If there was a cell $(r,c_1)$ with $r_1 < r < r_2$,
then the northeast property of $D_0$ would imply the presence of a cell in column $c_2$ that would contradict our choice of $x_2$. This leaves us to consider position $(r,c_1)$ in the following three cases: $r_1=r_2$ with $r=r_3$ or $r_1<r_2$ with $r =r_2,\ r_3$.

\textbf{Case 1:} $r_1=r_2$, $(r_3, c_1)\in D_0$. Let $x_4 := (r_3, c_1)$. If there were any cell $(r,c)$ with $r_1<r<r_3$ and $c<c_2$, by the northeast property we would have $(r,c_2)\in D_0$, contradicting our choice of $x_3$ as having minimal row index.  Thus, this rectangular region with opposite corners $(r_1+1,c_1)$ and $(r_3-1,c_2)$ is empty. 

Let $\tilde{x}_4=(\tilde{r}_4,c_1)$ and $x_4'=(r'_4,c_1)$ denote the positions of the $r_3$-cells in column $c_1$ of the standard labeling of $\tilde{D}$ and $D_{can}$,
respectively. Note that $r_{\max} \geq \tilde{r}_1>r_1'\geq r_2'$ since we have assumed $r_1=r_2$. Since column $c_2$ has rows $1$ through $r_{\max}$ occupied in $D_{can}$ (and $\tilde{D}$), it must be the case that $r_3'= r_2'+1$, so this implies that $r_3' \leq r_{\max}$.
By the emptiness of the rectangle described above, applying the canonical procedure to $D_0$ yields $r_4'=r_1'+1$.
Thus, if $\tilde{r}_1>r_1'$, we must have $\tilde{r}_4>r_4'$ as well. In this case, $\tilde{r}_4\leq r_{\max}$ means $\tilde{x}_4$ must 
have a cell $\tilde{y}$ in $\tilde{D}$ in column $c_2$ to its right, where $\tilde{y}$ would have started as cell $y$ above $x_3$ in $D_0$. Thus, we have shown that the triple $x_4$, $x_3$, $y$ satisfies (a) and (b) as well as the property $\tilde{r}_4>r_4'$, contradicting the maximality of our original triple choice.

\textbf{Case 2(a):} $r_1<r_2$, $(r_2, c_1)\in D_0$. If $D_0$ had a cell in column $c_2$ strictly below row $r_1$, then the northeast property of $D_0$ gives $(r_1,c_2)\in D_0$ which contradicts the maximality of the triple $x_1$, $x_2$, and $x_3$.  If $D_0$ had a cell in column $c_2$ strictly below row $r_2$ and weakly above row $r_1$, then again maximality is contradicted. Thus, column $c_2$ in $D_0$ has no cells below row $r_2$. Let $x_4=(r_2,c_1)\in D_0$.

In this case, instead of showing that the presence of $x_4$ leads to contradiction, we show that its presence can be used to deduce the existence of a different triple satisfying (a) - (e). Observe that combining $\cwt(D_0)_{c_2}>\cwt(D_0)_{c_1}$ and $\cwt_{\leq r_2}(D_0)_{c_1}\geq 2$, we have
\begin{equation*}
\cwt(D_0)_{c_2}\geq 3+\cwt_{>r_2}(D_0)_{c_1}.\label{eq:bdd-e-case-2a}\end{equation*}
But since $\cwt(D_0)_{c_2}=1+\cwt_{>r_2}(D_0)_{c_2}$, this equation becomes
\begin{equation}\label{eq:bdd-e-case-2a1}
\cwt_{>r_2}(D_0)_{c_2}\geq 2+\cwt_{>r_2}(D_0)_{c_1}.
\end{equation}

This implies that there is some row index $t_3 >r_2$ 
such that $z_3:=(t_3,c_2)\in D_0$ while $(t_3,c_1)\not\in D_0$. Let $t_1$ be the largest row index satisfying $r_2 \leq$ $t_1<t_3$ and $z_1:=(t_1,c_1)\in D_0$. Then by the northeast property of $D_0$, $z_2:=(t_1,c_2)\in D_0$ as well since $t_1 \geq r_2$.
Now, $z_1$, $z_2$, $z_3$ is a triple that satisfies (a), (b) and (e) by construction, as well as (c) and (d) since the columns have not changed from our original triple. 

\textbf{Case 2(b):} $r_1<r_2$, $(r_3, c_1)\in D_0$. As in the previous case, we know that column $c_2$ is empty below row $r_2$, since otherwise the minimality of our triple $x_1$, $x_2$, and $x_3$ is violated. We give a sketch of this case since it is nearly identical to the previous case.

Here, we instead have $\cwt_{\leq r_3}(D_0)_{c_1}\geq 2$ and $\cwt(D_0)_{c_2}=2+\cwt_{>r_3}(D_0)_{c_2}$, so the analog of \eqref{eq:bdd-e-case-2a1} becomes
\[\cwt_{>r_3}(D_0)_{c_2}\geq 1+\cwt_{>r_3}(D_0)_{c_1}.\]
This can similarly be used to deduce the existence of a cell $z_3\in D_0$ somewhere above row $r_3$, from which we obtain $z_1$ and $z_2$ in an identical way as in the previous case. The Properties (a) - (e) apply to the triple $z_1$, $z_2$, $z_3$ for the exact same reasons as the previous case. 

Having shown that all conditions (a) - (e) are satisfied by either the triple $x_1$, $x_2$, $x_3$ or a related triple in the same columns, this concludes the proof of Proposition~\ref{cor:mainNEboundedredux} and then Theorem~\ref{thm:mainNEbounded}.
\end{itemize}\end{proof}

\section{Applications to lock diagrams}\label{sec:lock}

We are now equipped to apply our main theorems to lock diagrams, a special class of northeast diagrams. First, we specialize Theorem \ref{thm:mainMFree} to lock polynomials. To do this, we give one additional definition.

\begin{definition}
A \textbf{subcomposition} of a weak composition $\alpha \in \comp_n$ is a weak composition $\beta \in \comp_k$ for $k \leq n$ such that there exist indices $1 \leq i_1 < i_2 < \cdots < i_k \leq n$ with $\beta_{j} = \alpha_{i_j}$ for all $1 \leq j \leq k$.
\end{definition}

\begin{coro}
    The lock polynomial $\mathfrak{K}_{\lock(\alpha)}$ is monomial multiplicity-free if and only if $\alpha$ has no subcomposition of the form $(0, 0, \alpha_i, \alpha_j)$ for $\alpha_i > 1$ and $\alpha_j > 0$.
\end{coro}

\begin{proof}
    (Proof of $\Rightarrow$) Suppose $\alpha$ has a subcomposition of the form $(0, 0, \alpha_i, \alpha_j)$ with $\alpha_i > 1$ and $\alpha_j > 0$. Then take $x_1$ to be the leftmost cell in row $i$ and $x_2$ the rightmost cell in row $j$ of $\lock(\alpha)$. We will show that $x_1$ and $x_2$ are the cells described in Theorem \ref{thm:mainMFree} by verifying each of the required conditions:
    \begin{itemize}
        \item[(a)] This is a result of the fact that $i < j$.
        \item[(b)] Because $\alpha_i > 1$, $x_1$ is not the rightmost cell in its row. Then $x_1$ is strictly to the left of $x_2$.
        \item[(c)] There are at least two empty rows below row $i$, so this condition is satisfied.
        \item[(d)] Again there are at least two empty rows below row $i$, so this condition is satisfied.
    \end{itemize}
    
    (Proof of $\Leftarrow$) Suppose $\mathfrak{K}_{\lock(\alpha)}$ is not monomial multiplicity-free. Then the cells $x_1$ and $x_2$ as described in Theorem \ref{thm:mainMFree} exist in $\lock(\alpha)$. Let $i$ be the row of $x_1$ and $j$ the row of $x_2$. Because of condition $(a)$, $i < j$. Because $x_2$ is a cell in row $j$, row $j$ is certainly nonempty, so $\alpha_j > 0$. Furthermore, because $x_1$ is to the left of $x_2$, $\alpha_i > 1$. Finally, because of condition $(d)$, there are at least two empty positions below the rightmost cell in row $i$ of $\lock(\alpha)$. As $\lock(\alpha)$ is a southeast diagram, the corresponding rows must be entirely empty. Then $\alpha$ has the desired subcomposition.
\end{proof}

Next we specialize Theorem \ref{thm:mainNEranked} to lock diagrams.

\begin{coro}\label{cor: ranked_locks}
    $\mathcal{P}(\lock(\alpha))$ is ranked if and only if for every pair $\alpha_i, \alpha_{i+k} \geq 2$ with $\alpha_{i+j} \in \{0,1\}$ for all $1 \leq j < k$, we have $\#\{j :  1 \leq j < k \text{ and } \alpha_{i+j}=1 \} \geq \#\{j :  j < i\text{ and } \alpha_{j}=0 \}$.
\end{coro}
\begin{proof}
As lock diagrams are a subset of northeast diagrams, we know that $\mathcal{P}(\lock(\alpha))$ is ranked if and only if it does not contain $x_1, x_2, x_3$ satisfying the conditions in Theorem \ref{thm:mainNEranked}. We show that in a lock diagram, having this forbidden configuration is equivalent to having a pair of row weights $\alpha_i, \alpha_{i+k} \geq 2$ with $\alpha_{i+j} \in \{0,1\}$ for all $1 \leq j < k$ such that $\#\{j :  1 \leq j < k \text{ and } \alpha_{i+j}=1 \} < \#\{j :  j < i\text{ and } \alpha_{j}=0 \}$. 

First assume that $\lock(\alpha)$ has $x_1, x_2,$ and $x_3$ satisfying conditions (a), (b), (c), and (d) of Theorem \ref{thm:mainNEranked}. In particular, consider such a configuration where $x_2$ and $x_3$ are chosen to have the minimum possible row index with respect to the position of $x_1$. As $\lock(\alpha)$ is right-justified, $x_2$ and $x_3$ lie in the same row. Set $x_1=(r_1,c_1)$, $x_2=(r_2,c_1)$, and $x_3=(r_2,c_3)$ as in the statement of the theorem. Then all positions to the right of cells $x_1$ and $x_2$ in rows $r_1$ and $r_2$, respectively, are occupied by cells.

If $c_n$ is the rightmost column of $\lock(\alpha)$, then we know $c_n$ satisfies condition (d), so there is at least one  entirely empty row below row $r_1$. Thus, all columns $c \geq c_1$ satisfy condition (c). Therefore, we can instead consider the triple $x_1=(r_1,c_{n-1})$, $x_2=(r_2,c_{n-1})$, $x_3=(r_2,c_n)$, which satisfies all of the conditions as well. If $r_2$ is not the first row of length at least 2 above row $r_1$, we can now reassign $r_2$ to be this row, and the conditions remain satisfied. 

With this final choice of $x_1$, $x_2$, and $x_3$, we see that 
condition (d) implies that the number of length 1 rows lying strictly between rows $r_1$ and $r_2$ is less than the number of empty rows below row $r_1$. Setting $r_1 = i$ and $r_2 = i+k$, we conclude that in a lock diagram, having the forbidden configuration described in Theorem \ref{thm:mainNEranked} implies that we have a pair of row weights $\alpha_i, \alpha_{i+k} \geq 2$ with $\alpha_{i+j} \in \{0,1\}$ for all $1 \leq j < k$ such that $\#\{j :  1 \leq j < k \text{ and } \alpha_{i+j}=1 \} < \#\{j :  j < i\text{ and } \alpha_{j}=0 \}$. 

For the converse, assume $\lock(\alpha)$ has rows $i$ and $i+k$ such that $\alpha_i, \alpha_{i+k} \geq 2$, $\alpha_{i+j} \in \{0,1\}$ for all $1 \leq j < k$, and $\#\{j :  1 \leq j < k \text{ and } \alpha_{i+j}=1 \} < \#\{j :  j < i\text{ and } \alpha_{j}=0 \}$. Let $c_n$ be the rightmost column of $\lock(\alpha)$ and choose $x_1 = (r_i, c_{n-1})$, $x_2=(r_{i+k},c_{n-1})$, and $x_3=(r_{i+k}, c_n)$. Then this triple of cells satisfies conditions (a) (b), (c), and (d) by construction.
\end{proof}

Specializing Theorem \ref{thm:mainNEbounded} to lock diagrams, we have the following result.

\begin{coro}\label{cor: bounded_locks} 
 $\mathcal{P}(\lock(\alpha))$ is bounded if and only if in the row weight vector $\alpha$, the nonzero entries after the first $0$ are weakly increasing.
\end{coro}
\begin{proof}
Every lock diagram is northeast, so $\mathcal{P}(\lock(\alpha))$ is bounded if and only if it does not contain $x_1, x_2, x_3$ satisfying the conditions (a) - (e) outlined in Theorem \ref{thm:mainNEbounded}. We show that in a lock diagram, having this forbidden configuration is equivalent to having some $i < j < k$ such that $\alpha_i =0$ and $\alpha_j > \alpha_k>0$.

As $\lock(\alpha)$ is right-justified, we observe that the set of rows occupied by cells in a column $c_k$ is always a subset of the set of rows occupied by cells in column $c_{k+1}$. Thus, condition (c) is true of any lock diagram where we choose columns $c_1 < c_2$ such that column $c_2 -1$ is not identical to column $c_2$. Thus, condition (c) holds exactly when there is some pair of indices $j,k$ with $\alpha_j\neq\alpha_k$. 

Furthermore, in a lock diagram, the existence of $x_1, x_2, x_3$ satisfying (a) and (b) is equivalent to the existence of an $\alpha_j>2$ and an $\alpha_k>1$ for some $k>j$, since we can always take $r_2=r_1=j$ and $r_3=k$. 
Adding condition (e) to this configuration, we see it is satisfied if and only if $(r_3,c_1)$ is empty, meaning that row $r_1$ has a greater weight than row $r_3$ or $\alpha_j>\alpha_k$.

Finally, condition (d) in a lock diagram is equivalent to saying that there is an empty row below row $r_1=j$. That is, there is $i$ with $i<j$ and $\alpha_i=0$. Thus, we conclude that for lock diagrams, the forbidden configuration for boundedness in northeast Kohnert posets is equivalent to having some $i < j < k$ such that $\alpha_i =0$ and $\alpha_j > \alpha_k>0$.
\end{proof}

We can now combine Corollaries \ref{cor: bounded_locks} and \ref{cor: ranked_locks} to give a classification of lock diagrams that yield ranked and bounded Kohnert posets.

\begin{coro} \label{cor: lock_rkd+bdd}
    $\mathcal{P}(\lock(\alpha))$ is ranked and bounded if and only if $\alpha$ satisfies the following: If $i$ is the smallest index such that $\alpha_i = 0$, then for all $j>i$, we have $\alpha_j>1$ if and only if $\alpha_j$ is the last nonzero entry in $\alpha$.     
\end{coro}

\begin{figure}[ht]
  \begin{center}
    \begin{tikzpicture}[scale=0.5]
    \cel(0,0) \cel(1,0) \cel(2,0)
    \cel(1,1) \cel(2,1)
    \cel(2,3) \cel(2,4) \cel(2,5)
    \cel(0,7) \cel(1,7) \cel(2,7)
    \draw (0,9)--(0,0)--(3,0);
    \end{tikzpicture}
    \caption{An example of a lock diagram that yields a ranked and bounded Kohnert poset.
    }\label{fig: lock_rkd+bdd}
  \end{center}
\end{figure}
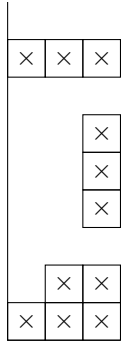

In other words, $\alpha = (\alpha_1,\ldots,\alpha_{\ell})$ has the following form:
\begin{itemize}
    \item $\alpha_1$ through $\alpha_{i-1}$ are greater than $0$, 
    \item $\alpha_i = 0$, 
    \item $\alpha_{i+1}$ through $\alpha_{\ell-1}$ are either $0$ or $1$, and
    \item $\alpha_{\ell} > 0$. 
\end{itemize}
An example of such a diagram can be found in Figure \ref{fig: lock_rkd+bdd}.

\begin{proof} For $\mathcal{P}(\lock(\alpha))$ to be bounded, recall that $\alpha$ must be such that the nonzero entries that follow the first $\alpha_i = 0$ are weakly increasing. Thus, $\alpha$ is such that $\alpha_1$ through $\alpha_{i-1}$ are greater than 0, $\alpha_i =0$, and for each $\alpha_j \neq 0$ with $j>i$ and $\alpha_k$ being the next nonzero entry following $\alpha_j$, we have $\alpha_j \leq \alpha_k$. 

Recall also that for $\mathcal{P}(\lock(\alpha))$ to be ranked,  we must have that for every pair $\alpha_j, \alpha_{j+k} \geq 2$, the number of 1s between them must be greater than or equal to the number of 0s below $\alpha_j$. Considering $\alpha$ for which $\mathcal{P}(\lock(\alpha))$ is bounded, we observe that this condition holds by default for $j < i$. For $j > i$, if $\alpha_j, \alpha_{j+k} \geq 2$, then since $\mathcal{P}(\lock(\alpha))$ is bounded, there cannot be any 1s between $\alpha_j$ and $\alpha_{j+k}$. But we know that $\alpha_i=0$ and $i<j$, so we conclude that if $\mathcal{P}(\lock(\alpha))$ is both ranked and bounded, then we cannot have two entries $\alpha_j, \alpha_{j+k} \geq 2$ if $j > i$.  

Thus, at most one nonzero entry after $\alpha_i$ is greater than 1, and by the weakly increasing condition for boundedness, this entry must be the final one, $\alpha_\ell$. We conclude that if $\mathcal{P}(\lock(\alpha))$ is both ranked and bounded, then it has the aforementioned form. 
\end{proof}

\section{Acknowledgments}

This project began during the 2024 Graduate Research Workshop in Combinatorics, which was supported by the University of Wisconsin-Milwaukee, the Combinatorics Foundation, and the National Science Foundation (NSF Grant DMS-1953445). B. A. Castellano was supported by an NSF Graduate Research Fellowship. A. Bingham was partially supported by supported by the grant DST/INT/RUS/RSF/P-41/2021 from the Department of Science \& Technology, Govt. of India, and FONDECYT-ANID grant 3250472. We would like to thank Kim Harry, Joakim Jakovleski, and Chelsea Sato for their early contributions to this project. We would also like to thank Nick Mayers for helpful discussions related to this work.

\printbibliography
\end{document}